\documentclass[11pt]{article}

\usepackage{amsmath,amsthm,amssymb,graphicx,xcolor}
\usepackage{fullpage}
\usepackage{mathtools}
\usepackage{color}
\usepackage{enumitem}
\usepackage{comment}
\usepackage{setspace}
\usepackage[colorlinks,citecolor=blue,urlcolor=blue]{hyperref}

\newtheorem{theorem}{Theorem}[section]
\newtheorem*{theorem*}{Theorem}
\newtheorem{proposition}[theorem]{Proposition}
\newtheorem*{proposition*}{Proposition}
\newtheorem{lemma}[theorem]{Lemma}
\newtheorem{claim}[theorem]{Claim}
\newtheorem{corollary}[theorem]{Corollary}
\newtheorem{definition}[theorem]{Definition}
\newtheorem{remark}[theorem]{Remark}
\newtheorem{assumption}[theorem]{Assumption}
\theoremstyle{remark}
\newtheorem*{remark*}{Remark}

\newcommand{\R}{\mathbb{R}}
\newcommand{\N}{\mathbb{N}}

\newcommand{\C}{\mathbb{C}}
\newcommand{\D}{\mathbb{D}}
\newcommand{\Dc}{\mathcal{D}}
\newcommand{\Ec}{\mathcal{E}}

\newcommand{\Qc}{\mathcal{Q}}

\newcommand{\Expect}[1]{\mathbb{E} \left[ #1 \right] }

\newcommand{\abs}[1]{\left\vert #1 \right\vert}
\newcommand{\norme}[1]{\left\| #1 \right\| }
\newcommand{\scalar}[1]{\left\langle #1 \right\rangle }
\newcommand{\floor}[1]{\left\lfloor #1 \right\rfloor}

\newcommand{\eps}{\varepsilon}

\newcommand{\E}{\mathbb{E}}

\renewcommand{\P}{\mathbb{P}}

\newcommand{\T}{\mathbb{T}}
\newcommand{\reals}{\mathbb{R}}
\newcommand{\complexes}{\mathbb{C}}

\newcommand{\integers}{\mathbb{Z}}

\newcommand{\1}{\mathbf{1}}
\newcommand{\HC}{H_{\complexes}}

\DeclareMathOperator{\supp}{supp}

\DeclareMathOperator{\Ker}{Ker}
\DeclareMathOperator{\Bin}{Bin}
\renewcommand{\Re}{\operatorname{Re}}
\renewcommand{\Im}{\operatorname{Im}}

\title{Density of imaginary multiplicative chaos via Malliavin calculus}
\author{Juhan Aru, Antoine Jego, Janne Junnila}

\date {}

\numberwithin{equation}{section}

\begin{document}

\renewcommand{\theparagraph}{\thesubsection.\arabic{paragraph}} 
\maketitle

\begin{abstract}
We consider the imaginary Gaussian multiplicative chaos, i.e. the complex Wick exponential $\mu_\beta := :e^{i\beta \Gamma(x)}:$ for a log-correlated Gaussian field $\Gamma$ in $d \geq 1$ dimensions. We prove a basic density result, showing that for any nonzero continuous test function $f$, the complex-valued random variable $\mu_\beta(f)$ has a smooth density w.r.t. the Lebesgue measure on $\C$. As a corollary, we deduce that the negative moments of imaginary chaos on the unit circle do not correspond to the analytic continuation of the Fyodorov-Bouchaud formula, even when well-defined.

Somewhat surprisingly, basic density results are not easy to prove for imaginary chaos and one of the main contributions of the article is introducing Malliavin calculus to the study of (complex) multiplicative chaos. To apply Malliavin calculus to imaginary chaos, we develop a new decomposition theorem for non-degenerate log-correlated fields via a small detour to operator theory, and obtain small ball probabilities for Sobolev norms of imaginary chaos.

\end{abstract}
\onehalfspacing

\section{Introduction}

In this paper we study imaginary Gaussian multiplicative chaos, formally written as $\mu_\beta := :e^{i\beta \Gamma(x)}:$, where $\Gamma$ is a log-correlated Gaussian field on a bounded domain $U \subset \reals^d$ and $\beta$ a real parameter. The study of imaginary chaos can be traced back to at least \cite{DERR,Big}, in case of cascade fields to \cite{BJM}, and to \cite{LCRV, JSW} in a wider setting of log-correlated fields. 

Imaginary multiplicative chaos distributions $:e^{i\beta \Gamma(x)}:$ can be rigorously defined as distributions in a Sobolev space of sufficiently negative index \cite{JSW}. In the case where $\Gamma$ is the 2D continuum Gaussian free field (GFF), they are related to the sine-Gordon model \cite{LRV, JSW} and the scaling limit of the spin-field of the critical XOR-Ising model is given by the real part of $:e^{i2^{-1/2}\Gamma(x)}:$ \cite{JSW}. Imaginary chaos has also played a role in the study of level sets of the GFF \cite{SSV}, giving a connection to SLE-curves.
In \cite{CGPR} it was shown using Wiener chaos methods that certain fields constructed using the Brownian Loop Soup converge to imaginary chaos. 
Recently, reconstruction theorems have been proved for both the continuum \cite{AruJunnila} and the discrete version \cite{GarbanSepulveda} of the imaginary chaos, showing that, somewhat surprisingly, when $d \ge 2$ it is possible to recover the underlying field from the information contained in the imaginary chaos in the whole subcritical phase $\beta \in (0,\sqrt{d})$. 

In a wider context, real multiplicative chaos $:e^{\gamma \Gamma(x)}:$, with $\gamma \in \R$ has been the subject of a lot of recent progress (see e.g. reviews \cite{RVreview, EPreview}). Complex and in particular imaginary multiplicative chaos appear then naturally, for example, as analytic extensions in $\gamma$. Complex variants of multiplicative chaos also come up when studying the statistics of zeros of the Riemann zeta function on the critical line \cite{SW}.

The main result of this paper is the existence and smoothness of density for random variables of the type $\mu_\beta(f)$. The main contribution, however, is probably the technique used to prove the main result. Indeed, whereas in the case of imaginary multiplicative cascades \cite{BM} and real multiplicative chaos \cite{RV} rather direct Fourier methods give the existence of a density, this approach is problematic in the case of imaginary chaos. The main obstacle is the presence of cancellations that are difficult to control without an exact recursive independence structure or monotonicity. We circumvent these problems by turning to Malliavin calculus. Interestingly, in order to apply methods of Malliavin calculus we have to first obtain new decomposition theorems for log-correlated fields, and prove quite technical concentration estimates for tails of imaginary chaos.

\subsection{The main result: existence of density}
Let us now denote by $\mu = \mu_\beta$ the imaginary chaos with parameter $\beta \in (0,\sqrt{d})$ in $d$ dimensions. In the appendix of \cite{LSZ} and in \cite{JSW} the tails of this random variable were studied and it was shown that $\P[|\mu(f)| > t]$ behaves roughly like $\exp(-t^{2d/\beta^2})$ -- this basically follows from the fact that using Onsager inequalities, one can obtain a very good control on the moments of imaginary chaos. 

In the present article we are interested in the local properties of the law of $\mu_\beta(f)$ and our main result is that this random variable has a smooth density. The following slightly informal statement is made precise in Theorem \ref{thm:main}.

\begin{theorem*}
Let $\Gamma$ be a non-degenerate log-correlated field in an open domain $U$ and let $f$ be a nonzero continuous function with compact support in $U$. Then the law of $\mu_\beta(f)$ is absolutely continuous with respect to the Lebesgue measure on $\complexes$ and the density is a Schwartz function.

Moreover, for any $\eta > 0$ the density is uniformly bounded from above for $\beta \in (\eta, \sqrt{d})$ and converges to zero pointwise as $\beta \to \sqrt{d}$.

Finally, the same holds in the case where $\mu_\beta$ is the imaginary chaos corresponding to the field $\hat{\Gamma}$ with covariance $\E[\hat{\Gamma}(x) \hat{\Gamma}(y)] = -\log |x-y|$ on the unit circle, with $f$ being any nonzero continuous function defined on the circle.
\end{theorem*}

\begin{remark*}
The reason why the circle field is brought out separately is because it does not satisfy our definition of non-degenerate log-correlated fields, see Section \ref{sec:basic}, and requires a bit of extra work. With similar work other cases of degenerate log-correlated fields could be handled. However, a unified approach to handle a more general class of log-correlated fields is still lacking.

The requirement of compact support for $f$ can also be dropped in many situations. For example, the theorem is also true in the case where $\Gamma$ is the zero-boundary GFF on a bounded simply connected domain in $\reals^2$ and $f \equiv 1$.
\end{remark*}

This theorem has already proved to be useful in further study of imaginary chaos\footnote{A work in preparation studies the monofractal structure of imaginary chaos.}, but we also expect this basic result and the method to be useful more generally in the study of complex chaos \cite{LCRV}, and in studying the integrability results related to multiplicative chaos \cite{FB, DOZZ} and the Sine-Gordon model. Not only should one be able to use this technique to prove density results in these more general cases, but as a corollary one can deduce the existence of certain negative moments, which have played important role in the above-mentioned results. In a follow-up work, we will prove by independent methods that the density for imaginary chaos is in fact everywhere positive.

\subsection{An application to the Fyodorov-Bouchaud formula}

Let us mention here one direct application of our results, linking our studies to recent integrability results on the Gaussian multiplicative chaos stemming from Liouville conformal field theory \cite{DOZZ, FB}. Namely, in \cite{FB} the author proved that for real $\gamma \in (0,\sqrt{2})$ the total mass of $:e^{\gamma \widehat \Gamma(x)}:$, where $\widehat \Gamma$ is the log-correlated Gaussian field on $S^1$ with covariance $C(x,y) = - \log |x-y|$, has an explicit density w.r.t. the Lebesgue measure; this was conjectured in \cite{FBoriginal} and proved by different methods in \cite{ChhNaj}. Moreover, in Theorem 1.1 of \cite{FB} the author proves an explicit expression for the $p-$th moment of $Y_\gamma := \frac{1}{2\pi}\int_{S^1} :e^{\gamma \widehat \Gamma(x)}: dx$ with $-\infty < p < 2/\gamma^2$: 
\begin{equation}\label{EQ:FB}
\E \left(Y_\gamma^p\right) = \frac{\Gamma(1-p\gamma^2/2)}{\Gamma(1-\gamma^2/2)^p},
\end{equation}
where with a slight abuse of notation $\Gamma$ is here the usual $\Gamma$-function.\footnote{Notice that in that paper the author is using a different normalization of the field with local behaviour of $-2\log |x-y|$.} Notice that for any $p$, the expression is analytic in $\gamma$ (outside of isolated singularities) and in particular analytic in a neighbourhood around the imaginary axis. So naively one might think that at least as long as the moments are defined for $:e^{i\beta \widehat \Gamma(x)}:$, they would correspond to the expression given by \eqref{EQ:FB} with $\gamma = i\beta$. And indeed, it is not hard to see that for $p \in \N$ this is the case. Our results however imply that this cannot be true in general, even in the case where the $p-$th moment is well-defined for the imaginary chaos. In other words, the analytic extension of the moment formulas is in general different from naively changing $\gamma$ in the Wick exponential.

\begin{corollary}
Let $\widehat \mu_\beta$ be the imaginary chaos corresponding to the log-correlated field $\widehat \Gamma$ on the unit circle. Then $\E\left( \widehat \mu_\beta(S^1)^{-1}\right)$ converges to zero as $\beta \to 1$. In particular, $\E \left(\widehat \mu_\beta(S^1)^{-1}\right)$ does not agree with the analytic continuation of Equation \eqref{EQ:FB} for $\gamma \in (-i, i)$. 
\end{corollary}

\begin{proof}
From Theorem \ref{thm:main} it follows that 
$$|\E\left( \widehat \mu_\beta(S^1)^{-1}\right)| \le \E\left( |\widehat \mu_\beta(S^1)|^{-1}\right) \to 0$$ 
as $\beta \to 1$. On the other hand a direct check shows that in Equation \eqref{EQ:FB}, the expression remains uniformly positive for $p = -1$, when we set $\gamma = i\beta$ and let $\beta \to 1$.
\end{proof}

\begin{remark}\label{rem:ref}
It might be ineteresting to take note that almost surely $Y_\gamma$ does have an analytic continuation in $\gamma$ to the unit disk of radius $\sqrt{2}$ around the origin. Moreovoer, from Theorem 1.1 in \cite{FB} we know that for $\gamma \in [0,2]$, the law of $Y_\gamma$ is equal to  $\frac{1}{\Gamma(1-\frac{1}{2}\gamma^2)}Y^{-\frac{\gamma^2}{2}}$, with $Y \sim Exp(1)$. One can then interpret the above corollary as saying that for $\gamma = i\beta$, the law of $Y_{i\beta}$ cannot be given by $\frac{1}{\Gamma(1+\frac{1}{2}\beta^2)}Y^{\frac{\beta^2}{2}}$, with $Y \sim Exp(1)$.
\end{remark}

\subsection{Other results: a decomposition of log-correlated fields and Sobolev norms of imaginary chaos}

As mentioned, our main tool in the proof of Theorem \ref{thm:main} is Malliavin calculus which is an infinite-dimensional differential calculus on the Wiener space introduced by Malliavin in the seventies \cite{Malliavin}. Whereas Malliavin calculus has been used to prove density results in various other settings \cite{Nualart}, we believe that it is a novel tool in the context of multiplicative chaos and could possibly have further interesting applications - e.g. in proving density results for more general models. In order to apply Malliavin calculus, we need to derive some results that could be of independent interest.

First, we derive a new decomposition theorem for non-degenerate log-correlated fields.
The following statement is more carefully formulated in Theorem \ref{thm:decomposition} and the proof has an operator-theoretic flavour.

\begin{theorem*}
  Let $\Gamma$ be a non-degenerate log-correlated Gaussian field on an open domain $U \subseteq \reals^d$ with covariance kernel given by $-\log |x-y| + g(x,y)$ and $g$ subject to some regularity conditions. Then, for every $V \Subset U$ we may write (possibly in a larger probability space)
  \[\Gamma|_V = Y + Z,\]
  where $Y$ is an almost $\star$-scale invariant field and $Z$ is a Hölder-regular field independent of $Y$, both defined on the whole of $\R^d$.
\end{theorem*}

Second, we develop a way to study the small ball probabilities of $\|f\mu_\beta\|_{H^{-d/2}(\reals^d)}$. The precise version of the following statement is given by Proposition \ref{prop:chaos_sobolev_tails}.

\begin{proposition*}
  Let $f \in C_c^\infty(U)$. Then for all $\beta \in (0, \sqrt{d})$ the probability $\P[\|f \mu_\beta\|_{H^{-d/2}(\reals^d)} \le  \lambda]$ decays super-polynomially as $\lambda \to 0$.
\end{proposition*}

This result is closely related to small ball probabilities of the Malliavin determinant of $\mu_\beta(f)$. To prove it we establish concentration results on the tail of imaginary chaos.

\subsection{Structure of the article}

We have set up the article to highlight how the general theory of Malliavin calculus is applied to prove such a density result and what are the concrete estimates of imaginary chaos needed to apply it. After collecting some preliminaries in Section \ref{sec:basic}, we use Section \ref{sec:malliavin} to walk the reader through the relevant notions and results of Malliavin calculus in the context of imaginary multiplicative chaos, thereby building up the backbone of the proof of the main theorem. In that section we state carefully the main result, and prove it up to technical estimates. The remaining proofs are then collected in Section \ref{sec:mvclc} and in Section \ref{sec:estimates_chaos}; the former contains some general lemmas of Malliavin calculus, and the latter deals with concentration results for imaginary chaos, including the proof of the Proposition \ref{prop:chaos_sobolev_tails} above. In Section \ref{sec:agd} we prove the decomposition theorem stated above. 

\subsection*{Acknowledgements}

We are thankful to two anonymous referees for their careful reading and their insights. In particular, Remark \ref{rem:ref} was added in the light of a referee's comments.
J.A.\ was supported by Eccellenza grant 194648 of the Swiss National Science Foundation and is a member of Swissmap. 
A.J. was recipient of a DOC Fellowship of the Austrian Academy of Sciences at the Faculty of Mathematics of the University of Vienna.
 
\section{Basic notions and definitions}\label{sec:basic}

\subsection{Log-correlated Gaussian fields and imaginary chaos}

\noindent In this section we establish the formal setup for the log-correlated field $\Gamma$ and of the imaginary chaos associated to $\Gamma$, often denoted by $:\exp(i\beta \Gamma):$ with $\beta \in \R$.

\subsubsection{Log-correlated Gaussian fields}

Let $U \subset \reals^d$ be a bounded and simply connected domain and suppose we are given a kernel of the form
\begin{equation}\label{eq:logkernel}
C(x,y) = \log \frac{1}{|x-y|} + g(x,y)
\end{equation}
where $g$ is bounded from above and satisfies $g(x,y) = g(y,x)$. Furthermore, we assume that $g \in H^{d+\varepsilon}_{\mathrm{loc}}(U \times U) \cap L^2(U \times U)$ for some $\eps > 0$.\footnote{For any $s \in \reals$ and $U \subset \reals^d$ we denote by $H^s_{\mathrm{loc}}(U)$ the space of distributions $f$ for which $\varphi f \in H^s(\reals^d)$ for all $\varphi \in C_c^\infty(U)$.} We may also extend $C(x,y)$ as $0$ outside of $U \times U$. Then $C$ defines a Hilbert--Schmidt operator on $L^2(\reals^d)$, and hence $C$ is self-adjoint and compact.

Assuming $C$ is positive definite, by spectral theorem there exists a sequence of strictly positive eigenvalues $\lambda_1 \ge \lambda_2 \ge \dots > 0$ and corresponding orthogonal eigenfunctions $(f_k)_{k \ge 1}$ spanning the subspace $L \coloneqq (\Ker C)^\bot$ in $L^2(\reals^d)$. We may now construct the log-correlated field $\Gamma$ with covariance kernel $C(x,y)$ via its Karhunen--Lo\`eve expansion
\begin{equation}\label{eq:logfieldexpansion}
\Gamma = \sum_{k \ge 1} A_k C^{1/2} f_k = \sum_{k \ge 1} A_k \sqrt{\lambda_k} f_k,
\end{equation}
where $(A_k)_{k \ge 1}$ is an i.i.d. sequence of standard normal random variables. It has been shown in \cite[Proposition~2.3]{JSW} that the above series converges in $H^{-\varepsilon}(\reals^d)$ for any fixed $\varepsilon > 0$.

From the KL-expansion one can see that heuristically $\Gamma$ is a standard Gaussian on the space $H_\Gamma \coloneqq C^{1/2} L$. The space $H:=H_\Gamma$ is called the Cameron--Martin space of $\Gamma$, and it becomes a Hilbert space by endowing it with the inner product $\langle f, g \rangle_H = \langle C^{-1/2} f, C^{-1/2} g \rangle_{L^2}$, where $C^{-1/2} f, C^{-1/2} g \in L$. This definition makes sense since $C^{1/2}$ is an injection on $L$. We will define the KL-basis $(e_k)_{k \ge 1}$ for $H$ by setting $e_k \coloneqq \sqrt{\lambda_k} f_k$, and we will also write $\langle \Gamma, h \rangle_H \coloneqq \sum_{k=1}^\infty A_k \langle h, e_k \rangle_H$ for $h \in H$. The left hand side in the latter definition is purely formal since $\Gamma \notin H$ almost surely.

Let us finally define what we mean by a non-degenerate log-correlated field in all of this paper. 

\begin{definition}[Non-degenerate log-correlated field]\label{def:logfield}
Consider a kernel $C_\Gamma(x,y) = C(x,y)$ from \eqref{eq:logkernel} and the associated log-correlated field $\Gamma$, given by \eqref{eq:logfieldexpansion}. We call the kernel $C$ and the field $\Gamma$ non-degenerate when $C$ is an injective operator on $L^2(U)$, i.e. $\Ker C = \{ 0 \}$. 
\end{definition}

Note that for covariance operators injectivity is equivalent to being strictly positive in the sense that $\langle C_\Gamma f, f \rangle > 0$ for all $f \in L^2(U)$, $f \neq 0$.\footnote{On $\reals^d$ one could also imagine a different definition of non-degenerate fields. Namely, a canonical way to define a log-correlated field $\Gamma_d$ on $\reals^d$ for any $d \ge 1$ is to take $H^{d/2}(\reals^d)$ as the Cameron--Martin space of the field. It would then be natural to call any log-correlated field on $\reals^d$ non-degenerate if its Cameron--Martin space is equivalent to $H^{d/2}(\reals^d)$. We will basically see in Section 4 that very roughly our condition implies that the Cameron--Martin space of a suitable extension of the non-degenerate field $\Gamma$ to the whole plane coincides up to an equivalent norm with $H^{d/2}(\reals^d)$.}\\

\noindent \textbf{The standard log-correlated field on the circle.}

The only degenerate field we will work with in this paper is the standard log-correlated field on the circle. I.e. it is the field $\Gamma$ on the unit circle which has the covariance $C_\Gamma(x,y) = \log \frac{1}{|x-y|}$, where one now thinks of $x$ and $y$ as being complex numbers of modulus $1$. 
Equivalently, we may consider the field on $[0,1]$ with the covariance
\[\E[\Gamma(e^{2\pi i t}) \Gamma(e^{2\pi i s})] = \log \frac{1}{2|\sin(\pi(t-s))|},\]
in which case we may write
\[\Gamma(e^{2 \pi i t}) = \sqrt{2} \sum_{k=1}^\infty \frac{1}{\sqrt{k}}(A_k \cos(2\pi k t) + B_k \sin(2\pi k t))\]
where $A_k$ and $B_k$ are i.i.d. standard normal random variables.

This circle field is degenerate because it is conditioned to satisfy $\int_0^1 \Gamma(e^{2 \pi i \theta}) \, d\theta = 0$ and the operator $C$ maps constant functions to zero. It is however not hard to see that after restricting the domain of the field $\Gamma(e^{2\pi i \cdot})$ to $I_0 \coloneqq [-1/4,1/4]$ it becomes non-degenerate.

\subsubsection{Imaginary chaos}

Let us now fix $\beta \in (0,\sqrt{d})$. For any $f \in L^\infty(U)$ we may define the imaginary chaos $\mu$ tested against $f$ via the regularization and renormalisation procedure
\[\mu(f) \coloneqq \lim_{\varepsilon \to 0} \int_U f(x) e^{i\beta \Gamma_\varepsilon(x) + \frac{\beta^2}{2} \E \Gamma_\varepsilon(x)^2} \, dx,\]
where $\Gamma_\varepsilon$ is a convolution approximation of $\Gamma$ against some smooth mollifier $\varphi_\varepsilon$. An easy computation shows that the convergence takes place in $L^2(\Omega)$. Importantly, the limiting random variable does not depend on the choice of mollifier. Again, one has to be careful however when defining $\mu(f)$ for uncountably many $f$ simultaneously. Indeed, $\mu$ turns out to have a.s. infinite total variation, but it does define a random $H^s(\reals^d)$-valued distribution when $s < -\beta^2/2$ \cite{JSW}. One may also (via a change of the base measure in the proofs of \cite{JSW}) fix $f \in L^\infty(\reals^d)$ and consider $g \mapsto \mu(fg)$ as an element of $H^s(\reals^d)$. Although $\mu$ is not defined pointwise, we will below freely use the notation $\int_U f(x) \mu(x) \, dx$ to refer to $\mu(f)$. 

\subsection{Malliavin calculus: basic definitions}

In this subsection we will collect some very basic notions of Malliavin calculus: the Malliavin derivative and Malliavin smoothness. We will mainly follow \cite{Nualart} in our definitions, making some straightforward adaptations for complex-valued random variables both here and in the following sections.
 
Let $C_p^\infty(\reals^n;\reals)$ be the class of real-valued smooth functions defined on $\reals^n$ such that $f$ and all its partial derivatives grow at most polynomially.
\begin{definition}\label{def:der}
  We say that $F$ is a smooth (real) random variable if it is of the form
\[F(\Gamma) = f(\langle \Gamma, h_1 \rangle_H, \dots, \langle \Gamma, h_n \rangle_H)\]
for some $h_1,\dots,h_n \in H$ and $f \in C_p^\infty(\reals^n;\reals)$, $n \ge 1$. 

For such a variable $F$ we define its Malliavin derivative $D F$ by
\[D F = \sum_{k=1}^n \frac{\partial}{\partial_k} f(\langle \Gamma, h_1 \rangle_H, \dots, \langle \Gamma, h_n \rangle_H) h_k.\]
\end{definition}
Thus we see that $DF$ is an $H$-valued random variable and in fact, in the case where $F$ is a smooth random variable, $DF$ corresponds to the usual derivative map: for any $h \in H$, we have that 
\[ \langle DF(\Gamma), h \rangle_H = \lim_{\eps \to 0} \frac{F(\Gamma + \eps h) - F(\Gamma)}{\eps}. \]
One may also define $D^mF$ as a $H^{\otimes m}$-valued random variable by setting
\[D^m F = \sum_{k_1,\dots,k_m=1}^n \frac{\partial^m}{\partial_{k_1} \dots \partial_{k_m}} f(\langle \Gamma, h_1 \rangle_H, \dots, \langle \Gamma, h_n \rangle_H) h_{k_1} \otimes \dots \otimes h_{k_m}.\]
In our case $H$ is a space of functions defined on $U$ and hence $H^{\otimes m}$ can be seen as a space of functions defined on $U^m$. At times it will be convenient to write down the arguments of the function explicitly using subscripts, e.g. for all $t_1,\dots,t_m \in U$ we set 
\[D^m_{t_1,\dots,t_m} F \coloneqq D^m F(t_1,\dots,t_m),\] with
\[D^m F(t_1,\dots,t_m) = \sum_{k_1,\dots,k_m=1}^n \frac{\partial^m}{\partial k_1 \dots \partial k_m} f(\langle \Gamma, h_1 \rangle_H, \dots, \langle \Gamma, h_n \rangle_H) h_{k_1}(t_1)\dots h_{k_m}(t_m).\]

We extend the above definition in a natural way to complex smooth random variables by setting
\[D (F + iG) = DF + i DG\]
when $F$ and $G$ are real smooth random variables. Thus in general $D$ will map complex random variables to the complexification of $H$, which we denote by $\HC$. We will assume that the inner product $\langle \cdot, \cdot \rangle_{\HC}$ is conjugate linear in the second variable. From here onwards we will use $F$ for complex-valued Malliavin smooth random variables, unless otherwise stated.

To define $D$ for a larger class of random variables one uses approximation by the smooth functions above. 
More precisely, we define for any non-negative integer $k$ and real $p \ge 1$ the class of random variables $\D^{k,p}$ as the completion of (complex) smooth random variables with respect to the norm
\[\|F\|_{k,p}^p \coloneqq \E |F|^p + \sum_{j=1}^k \E \|D^j F\|_{\HC^{\otimes j}}^p.\]
The spaces $\D^{k,p}$ are decreasing with $p$ and $k$, and we denote $\D^\infty \coloneqq \bigcap_{p,k \ge 1} \D^{k,p}$.
Similarly we set $\D^{k,\infty} \coloneqq \bigcap_{p \ge 1} \D^{k,p}$.

Finally, viewing $D$ as an unbounded operator on $L^2(\Omega;\complexes)$ with values in $L^2(\Omega;H_{\complexes})$, we may define its adjoint $\delta$ which is also called the divergence operator. More specifically we have 
\[\E[F \delta u] = \E \langle DF, u \rangle_{\HC}\]
for any $u$ such that $|\E \langle DF, u \rangle_{\HC}|^2 \lesssim \E F^2$ for all $F \in \D^{1,2}$.

\section{Density of imaginary chaos via Malliavin calculus}\label{sec:malliavin}

Let $f$ be a continuous function of compact support in $U$. Our goal is to apply Malliavin calculus to show that the random variable $M \coloneqq \mu(f)$ has a smooth density with respect to the Lebesgue measure on $\complexes$. 

We start by walking through the basic results of Malliavin calculus that we want to apply and we then reduce the proof of Theorem \ref{thm:main} to concrete estimates on imaginary chaos. Some useful lemmas of Malliavin calculus are proven in Section \ref{sec:mvclc} and the estimates on imaginary chaos are verified in Section \ref{sec:estimates_chaos}, with input from Section \ref{sec:agd}.

\medskip

Formally one can write the Malliavin derivative $DM$ of $M = \mu(f)$ as
\begin{align*}
  D_t M & = \int f(x) D_t :e^{i\beta \sum_{n=1}^\infty \langle \Gamma, e_n \rangle_H e_n(x)}: \, dx \\
            & = \int f(x) \sum_{k=1}^\infty :e^{i\beta \Gamma(x)}: i\beta e_k(t) e_k(x) \, dx \\
            & = i\beta \int f(x) \mu(x) C(t,x) \, dx. 
\end{align*}

The content of the following proposition is to make the above computations rigorous by truncating the series $\sum_{n=1}^\infty \langle \Gamma, e_n \rangle_H e_n(x)$ to be able to work with Malliavin smooth random variables, as in Definition \ref{def:der}. 
\begin{proposition}\label{prop:dinfty}
  Let $f \in L^\infty(\complexes)$. Then $M \in \D^{\infty}$ and
  \[D_t M = i\beta \int_U f(x) \mu(x) C(t,x) \, dx\]
  for all $t \in U$.
\end{proposition}

The reason we are interested in showing that $M$ belongs to $\D^\infty$ is the following classical result of Malliavin calculus, stating sufficient conditions for the existence of a smooth density.
For convenience we state it here directly for complex valued random variables.

\begin{proposition}
\label{prop:existence_density}
  Let $F \in \D^\infty$ be a complex valued random variable and let
  \begin{equation}
    \label{eq:malliavin_det}
    \det \gamma_F \coloneqq \frac{1}{4}(\|DF\|_{\HC}^4 - |\langle DF, D\overline{F}\rangle_{\HC}|^2)
  \end{equation}
  be the Malliavin determinant of $F$.
  If $\E |\det(\gamma_F)|^{-p} < \infty$ for all $p \ge 1$, then $F$ has a density $\rho$ w.r.t. the Lebesgue measure in $\complexes$ and $\rho$ is a Schwartz function.
\end{proposition}

The proof follows rather directly from  \cite[Proposition~2.1.5]{Nualart}:
\begin{proof}
  Following \cite{Nualart}, the Malliavin matrix of a random vector $F = (F_1,\dots,F_n) \in \reals^n$ is given by $\gamma_F \coloneqq (\langle DF_j,DF_k\rangle_H)^{n}_{j,k}$. We will use Proposition 2.1.5 from \cite{Nualart}, which states that if $F_i \in \D^\infty$ and $\E |\det \gamma_F|^{-p} < \infty$ for all $p \ge 1$, then $F$ has a density w.r.t. the Lebesgue measure on $\reals^n$ which is a Schwartz function.   
  
  As $\Re F, \Im F \in \D^\infty$ by assumption, it is enough to check that $\det \gamma_F$ is equal to the given formula in the case $F = (\Re F, \Im F)$. This is easy to check by writing
  \begin{align*}
    \det \gamma_F & = \langle DF_1, DF_1 \rangle_H \langle DF_2, DF_2 \rangle_H - \langle DF_1, DF_2, \rangle_H^2 \\
                  & = \frac{1}{16} \|DF + D\overline{F}\|_{\HC}^2 \|DF - D\overline{F}\|_{\HC}^2 - \frac{1}{16} |\langle DF + D\overline{F}, DF - D\overline{F}\rangle_{\HC}|^2
  \end{align*}
  and expanding the squares on the right hand side. We leave the details to the reader.
\end{proof}

Thus to show that $F$ has a smooth and bounded density it will be enough to show that the negative moments of $\|DF\|_{\HC}^4 - |\langle DF, D\overline{F}\rangle_{\HC}|^2$ are all finite. In fact this quantity is not straightforward to control directly and to make calculations possible, we first apply the following projection bounds, whose proofs we postpone to Section \ref{sec:mvclc}:

\begin{lemma}[Projection bounds]\label{lem:deltalowerbound}
	Let $F \in \D^{1,2}$ and let $h$ be any function in $\HC$. Then
	\begin{equation}\label{eq:deltalowerbound1}
	\frac{\det \gamma_F}{\|DF\|^2_{\HC}} \geq \frac{1}{4} \frac{(|\langle DF, h \rangle_{\HC}| - |\langle D\overline{F}, h \rangle_{\HC}|)^2}{\|h\|_{\HC}^2}.
	\end{equation}
	and
	\begin{equation}\label{eq:deltalowerbound2}
	\det \gamma_F \geq \frac{1}{4} \frac{(|\langle DF, h \rangle_{\HC}| - |\langle D\overline{F}, h \rangle_{\HC}|)^4}{\|h\|_{\HC}^4}.
	\end{equation}
	
\end{lemma}

To further show that the density is uniformly bounded in $\beta$ outside any interval surrounding the origin, we need to have some quantitative control on the densities. We will use the following simple adaption of Lemma 7.3.2 in \cite{NuaNua} to the complex case to do this: 
\begin{lemma}\label{lem:densityupbound}
Let $p > 2$ and $F$ be a complex Malliavin random variable in $\D^{2,\infty}$. Then there is a constant $c = c_p > 0$ depending only on $p$ such that the density $\rho$ of $F$ satisfies for all $x \in \C$
\[\rho(x) \leq c_p (\E|\delta(A)|^p)^{2/p},\]
where $A$ is defined by 
  \[A = \frac{\|DF\|_{\HC}^2 DF - \langle DF, D\overline{F} \rangle_{\HC} D\overline{F}}{\|DF\|_{\HC}^4 - |\langle DF, D\overline{F}\rangle_{\HC}|^2}.\]
\end{lemma}

Bounding $\delta(A)$ is again technically not straightforward, but the following general bound could possibly be of independent interest. It is again proved in Section \ref{sec:mvclc}.

\begin{proposition}\label{prop:bndcv}Let $F$ be a complex Malliavin random variable in $\D^{2,\infty}$. We have
  \[|\delta(A)| \lesssim \frac{\|DF\|_{\HC}^2(|\delta(DF)| + \|D^2F\|_{H_{\complexes} \otimes H_{\complexes}})}{\|DF\|_{\HC}^4 - |\langle DF, D\overline{F}\rangle_{\HC}|^2}.\]
\end{proposition}

Using the above results on Malliavin calculus, we can now reduce Theorem \ref{thm:main} to concrete propositions on imaginary chaos. Proving the estimates needed for these propositions is basically the content of Section \ref{sec:estimates_chaos}. 

\medskip 

We start with a precise statement of the main theorem:

\begin{theorem}\label{thm:main}
Let $U$ be an open bounded domain and $\Gamma$ a non-degenerate log-correlated field in $U$ as in Definition \ref{def:logfield} and $f$ be a nonzero continuous function of compact support in $U$. We denote by $\mu$ the imaginary chaos associated to $\Gamma$ and parameter $\beta \in (0,\sqrt{d})$. Then 
\begin{itemize}
    \item the law of $\mu(f)$ is absolutely continuous with respect to the Lebesgue measure on $\complexes$ and the density is a Schwartz function;
    \item for any $\eta > 0$ the density is uniformly bounded from above for $\beta \in (\eta, \sqrt{d})$ and converges to zero pointwise as $\beta \to \sqrt{d}$.
\end{itemize}
Finally, the same holds in the case where $\Gamma$ is defined on the unit circle with covariance $\E[\hat{\Gamma}(x) \hat{\Gamma}(y)] = -\log |x-y|$ and $f$ is any nonzero continuous function on the circle. 
\end{theorem}

There are basically two technical chaos estimates needed to deduce the theorem. First, super-polynomial bounds on small ball probabilities of the Malliavin determinant are used both to prove that the density exists and is a Schwartz function, and to show uniformity:

\begin{proposition}\label{prop:tail_det} Let $\Gamma$, $f$, $M = \mu(f)$ be as in the theorem above. Then we have the following bounds for the Malliavin determinant $\det \gamma_M$. For any $\nu > 0$, there exist constants $C, c, a, \eps_0 > 0$ (which do not depend on $\beta$) such that for all $\eps \in (0,\eps_0)$ and for all $\beta \in (\nu, \sqrt{d})$, 
\begin{equation}
\label{eq:prop_tail_det1}
\P \left( \det \gamma_M \geq (d-\beta^2)^{-4}\eps \right) \geq  1 - C\exp\left(-a \varepsilon^{-c/2}\right).
\end{equation}
and
\begin{equation}
\label{eq:prop_tail_det2}
\P \left( \frac{\det \gamma_M}{\|DM\|^2_{\HC}} \geq (d-\beta^2)^{-2}\eps \right) \geq  1 - C\exp\left(-a \varepsilon^{-c}\right).
\end{equation}
\end{proposition}

Here the bound on $\frac{\norme{DM}_{H_\C}^2}{\det \gamma_M}$ is necessary, when bounding the divergence of the covering field via Proposition \ref{prop:bndcv}. Second, in order to apply Lemma \ref{lem:densityupbound} we also need upper bounds on $|\delta(DM)|$ and $\|D^2M\|_{\HC \otimes \HC}$: 

\begin{proposition}\label{prop:bound_delta_D2_intro}
 Let $\Gamma$, $f$, $M = \mu(f)$ be as in the theorem above. Then for all $N \geq 1$, there exists $C = C(N)>0$ such that for all $\beta \in (0, \sqrt{d})$ 
\begin{equation}\label{eq:prop_bound_delta_DM}
    \E \left[ \abs{ \delta(DM) }^{2N} \right] \leq C (d-\beta^2)^{-3N}
\end{equation}
and
\begin{equation}\label{eq:prop_bound_D2M}
    \E \left[ \|D^2M\|_{\HC \otimes \HC}^{2N} \right] \leq C (d-\beta^2)^{-3N}.
\end{equation}
\end{proposition}

We can now prove Theorem \ref{thm:main} modulo these propositions. 

\begin{proof}[Proof of Theorem \ref{thm:main}]
To apply Proposition \ref{prop:existence_density} to prove that $M = \mu(f)$ has a density w.r.t. Lebesgue measure, and that moreover this density is a Schwartz function, we need to verify two conditions:
\begin{itemize}
    \item That $M \in \D^\infty$ -- this is the content of Proposition \ref{prop:dinfty};
    \item And that $\E|\det(\gamma_M)|^{-p} < \infty$ for all $p \geq 1$ -- this follows directly from the bound \eqref{eq:prop_tail_det1} in Proposition \ref{prop:tail_det}.
\end{itemize}

Finally, it remains to argue that the density is uniformly bounded from above for $\beta \in (\eta, \sqrt{d})$ for some fixed $\eta > 0$, and converges to zero pointwise on $\R^d$ as $\beta \to \sqrt{d}$. This follows from Lemma \ref{lem:densityupbound}, once we show that
$\E |\delta(A)|^4$ is uniformly bounded in $\beta \in (\eta, \sqrt{d})$ and tends to zero as $\beta \to \sqrt{d}$.
By Proposition \ref{prop:bndcv}
\[\E |\delta(A)|^4 \lesssim \E\Big|\frac{\|DM\|_{\HC}^2(|\delta(DM)| + \|D^2M\|_{H_{\complexes} \otimes H_{\complexes}})}{\|DM\|_{\HC}^4 - |\langle DM, D\overline{M}\rangle_{\HC}|^2}\Big|^4.\]
By using the inequality $(x+y)^4\lesssim x^4 + y^4$ and then Cauchy--Schwarz we have that
\[\E |\delta(A)|^4 \lesssim \sqrt{\E \Big|\frac{\|DM\|_{\HC}^2}{\det \gamma_M}\Big|^8\E|\delta(DM)|^8} + \sqrt{\E \Big|\frac{\|DM\|_{\HC}^2}{\det \gamma_M}\Big|^8\E|\|D^2M\|_{H_{\complexes} \otimes H_{\complexes}}|^8}.\]
We thus conclude from \eqref{eq:prop_tail_det2} in Proposition \ref{prop:tail_det} and Proposition \ref{prop:bound_delta_D2_intro}.
\end{proof}

\noindent The proofs of the above-mentioned chaos estimates appear in Section \ref{sec:estimates_chaos}. More precisely, 
\begin{itemize}
\item In Section \ref{sec:lemma:dinfty} we prove that $M$ is in $\D^\infty$, i.e. Proposition \ref{prop:dinfty}. This boils down to bounding moments of $DM$ and is a rather standard calculation. Similar computations with small improvements on existing estimates allow to prove Proposition \ref{prop:bound_delta_D2_intro} in Section \ref{sec:bounds_for_derivatives}.
\item In Section \ref{sec:mainproof}, we prove Proposition \ref{prop:tail_det}, which requires a novel approach. It is also in this subsection where we make use of the almost global decomposition theorem for non-degenerate log-correlated fields, proved in Section \ref{sec:agd}.
\end{itemize}
The missing general results of Malliavin calculus are proved in Section \ref{sec:mvclc}.

\section{Almost global decompositions of non-degenerate log-correlated fields}\label{sec:agd}

It is often useful to try to decompose the log-correlated Gaussian field $\Gamma$ on the open set $U \subset \reals^d$ as a sum of two independent fields $Y$ and $Z$, where $Y$ is in some sense canonical and easy to calculate with, and $Z$ is regular.
In \cite{JSW2} it was shown that such decompositions exist around every point $x_0 \in U$ when $g \in H_{\mathrm{loc}}^{s}(U \times U)$ for some $s > d$ and $Y$ is taken to be a so-called almost $\star$-scale invariant field.

Our goal in this section is to establish a more general variant of this decomposition theorem which removes the need to restrict to small balls and works in any subdomain $V \Subset U$ (we write $A \Subset B$ to indicate that $\overline{A} \subset B$) by simply assuming that $\Gamma$ is non-degenerate on $V$, meaning that $C_\Gamma$ defines an injective integral operator on $L^2(V)$, as explained in Section \ref{sec:basic}. 

In the context of the present article, the usefulness of this result is strongly interlinked with the following standard comparison result for Cameron--Martin spaces. In the case of Reproducing Kernel Hilbert spaces, this can be found for example in \cite{Aronz}.

\begin{lemma}\label{lem:CMdecomp}
  Let $Y$ and $Z$ be two independent distribution-valued Gaussian fields and denote $\Gamma = Y + Z$.
  Let $(H_\Gamma, \| \cdot \|_{H_\Gamma})$ and $(H_Y, \| \cdot \|_{H_Y})$ be the Cameron--Martin spaces of $\Gamma$ and $Y$ respectively.
  Then $H_Y \subset H_\Gamma$ and moreover for every $h \in H_{Y}$, we have that $\|h\|_{H_Y} \geq \|h\|_{H_\Gamma}$.
\end{lemma}

Basically, via this Lemma our decomposition allows to meaningfully transfer calculations on the initial field $\Gamma$ to easier ones on the almost $\star$-scale invariant fields $Y$, where Fourier methods become available. 

We will start by recalling the basic definitions related to $\star$-scale invariant and almost $\star$-scale invariant log-correlated fields. We then state the theorem and discuss heuristics, and finally prove the theorem in two last subsections. In this section all function spaces are the standard function spaces for real-valued functions, i.e. we don't need to consider their complexified counterparts.

\subsection{Overview of $\star$-scale and almost $\star$-scale invariant log-correlated fields}\label{sec:starscale}

To define $\star$-scale invariant and almost $\star$-scale invariant fields, we first need to pick a seed covariance $k$. For simplicity we will in what follows make the following assumptions on $k$:
\begin{assumption}\label{assumption:seed_covariance}
	The seed covariance $k \colon \reals^d \to \reals$ satisfies the following properties:
	\begin{itemize}
		\item $k(x) \ge 0$ for all $x \in \reals^d$ and $k(0) = 1$;
		\item $k(x) = k((|x|, 0, \dots, 0)) \eqqcolon k(|x|)$ is rotationally symmetric and  $\supp k \subset B(0,1)$,
		\item There exists $s > \frac{d+1}{2}$ such that $0 \le \hat{k}(\xi) \lesssim (1 + |\xi|^2)^{-s}$ for all $\xi \in \reals^d$.
	\end{itemize}
\end{assumption}

The fact that $k$ is supported in $B(0,1)$ yields the useful property that distant regions of the associated Gaussian field will be independent.

Let us also remark that an easy way to construct a seed covariance $k$ satisfying the above assumptions is to take a smooth, non-negative and rotationally symmetric function $\varphi$ supported in $B(0,1/2)$ with $\|\varphi\|_{L^2} = 1$ and then letting $k = \varphi * \varphi$ be the convolution of $\varphi$ with itself.

\begin{definition}
  Let $k \colon \reals^d \to \reals$ be as above.  The $\star$-scale invariant covariance kernel $C_X$ associated to $k$ is given by
  \[C_X(x,y) \coloneqq \int_0^\infty k(e^{u}(x-y)) \, du.\] 
  Similarly, the related \emph{almost} $\star$-scale invariant covariance kernel $C_Y = C_{Y^{(\alpha)}}$ associated to $k$ and a parameter $\alpha > 0$ is given by
  \[C_Y(x,y) \coloneqq \int_0^\infty k(e^{u}(x-y)) (1 - e^{-\alpha u}) \, du.\]
\end{definition}

We often use approximations $Y_\delta$ of $Y$, which can be defined via the stochastic integrals
\begin{equation}\label{eq:YWN}
Y_\delta(x) = \int_{\reals^d \times [0, \log \frac{1}{\delta}]} e^{du/2} \tilde{k}(e^u(t-x)) \sqrt{1 - e^{-\alpha u}} dW(t,u),
\end{equation}
where $W$ is the standard white noise on $\reals^{d+1}$ and
$\tilde{k}(x) = \mathcal{F}^{-1}{\sqrt{\mathcal{F} k}}(x)$ with $\mathcal{F}$ denoting the Fourier transform.

We also define the tail field $\hat{Y}_\delta \coloneqq Y - Y_\delta$, which decorrelates at distances bigger than $\delta$.
The following lemma then gives basic estimates on the covariance of this tail field. See Appendix \ref{appendixA} for the proof.

\begin{lemma}\label{lemma:star_scale_cov_2}
  There exists a constant $C > 0$ such that
  \[\E[\hat{Y}_\delta(x) \hat{Y}_\delta(y)] \leq\frac{\delta}{|x-y|}\]
  and 
 \[\E[\hat{Y}_\delta(x) \hat{Y}_\delta(y)] \geq \frac{\delta}{|x-y|} - C.\]
  Moreover $\E[\hat{Y}_\delta(x) \hat{Y}_\delta(y)] = 0$ whenever $|x-y| \geq \delta$.
\end{lemma}

\subsection{Statement of the theorem and the high level argument}

The main theorem of this section can be stated as follows.

\begin{theorem}\label{thm:decomposition}
  Let $\Gamma$ be a non-degenerate log-correlated Gaussian field on an open domain $U \subseteq \reals^d$ as in Definition \ref{def:logfield}. Assume further that the covariance kernel given by \eqref{eq:logkernel} satisfies $g \in H^{s}_{\mathrm{loc}}(U \times U)$ for some $s > d$.
  
  Then for every seed kernel $k$ satisfying Assumption~\ref{assumption:seed_covariance} and every $V \Subset U$, there exists $\alpha > 0$ (possibly depending on $V$) such that we may write (possibly in a larger probability space)
  \[\Gamma|_V = Y + Z,\]
  where $Y$ is an almost $\star$-scale invariant field with seed covariance $k$ and parameter $\alpha$ and $Z$ is a Hölder-regular field independent of $Y$, both defined on the whole of $\R^d$. Moreover, there exists $\varepsilon > 0$ such that the operator $C_Z$ maps $H^s(\reals^d) \to H^{s + d + \varepsilon}(\reals^d)$ for all $s \in [-d,0]$.
\end{theorem}

Notice that the 2D zero boundary Gaussian free field is a non-degenerate log-correlated field in the open disk. However, there is no hope to decompose it using an almost $\star$-scale invariant field on the whole of $\D$, so in that sense the above theorem is as global as you could hope.\footnote{This can be checked e.g. by considering the equality $C_\Gamma(x,y) = C_Y(x,y) + C_Z(x,y)$ at two points $x$ and $y$ with $y$ tending first to a fixed boundary point $z$ and then $x$ tending to the same point. In the limit one formally obtains $0 = \infty$.}

\begin{remark}
In \cite[Theorem~B]{JSW2} it was shown that even for a degenerate log-correlated field $\Gamma$, one can for any $x \in U$ find a ball $B(x,r(x))$, restricted to which $\Gamma$ is non-degenerate and can be decomposed as an independent sum of an almost star-scale invariant field and a Hölder-regular field. In this sense one can see Theorem \ref{thm:decomposition} as a generalization in the special case of non-degenerate fields. 
\end{remark}

Before going to the proof of Theorem~\ref{thm:decomposition}, let us try to illustrate the high level argument in terms of the following toy problem on the unit circle $\T = \{z \in \complexes : |z| = 1\}$: Let $\Gamma$ be a non-degenerate log-correlated field on $\T$ with covariance of the form $\log \frac{1}{|x-y|} + g(|x-y|)$, where now also the $g$ term only depends on the distance between the two points.
This means that we can write the covariance using the Fourier series
\[C_\Gamma(x,y) = \frac{g_0}{2} + \Re\sum_{n=1}^\infty (\frac{1}{n} + g_n) x^n y^{-n},\]
where
\[g_n \coloneqq \frac{1}{\pi} \int_{\T} g(|1 - x|) x^{-n} |dx|,\]
with $|dx|$ denoting the arc-length measure. As $\Gamma$ is assumed to be non-degenerate, we know that $\frac{1}{n} + g_n > 0$ for all $n \geq 1$.

The almost $\star$-scale field would correspond to a field with covariance of the form
\[C_{Y}(x,y) = \Re \sum_{n=1}^\infty (\frac{1}{n} - \frac{1}{n^{1 + \alpha}}) x^n y^{-n},\]
and thus the difference between the tail and the two covariances would be 
\[C_\Gamma(x,y) - C_{Y}(x,y) = \frac{g_0}{2} + \Re \sum_{n=1}^\infty (\frac{1}{n^{1 + \alpha}} + g_n) x^n y^{-n}.\]
It is now easy to see that if $g_n = O(n^{-s})$ for some $s > 1 + \alpha$, the coefficients in the above difference are positive for all large enough $n$. By further reducing $\alpha$, we can guarantee that $\frac{1}{n^{1+\alpha}} + g_n > 0$ for all $n \geq 1$, so that the difference $C_\Gamma - C_Y$ is again a positive definite kernel.

The main issue in implementing this strategy for general log-correlated covariances on domains in $\reals^d$ is the fact that in general we do not have a canonical basis such that $C_\Gamma$ and $C_X$ would be simultaneously diagonalizable. To still be able to make useful calculations, we thus want to find some universal, non-basis dependent setting, where both can be studied. This is comfortably offered for example by the Fourier transform on spaces $L^2(\reals^d)$ and $H^s(\reals^d)$. Thus as a first step we will find a suitable extension of $\Gamma$ to a log-correlated field on the whole of $\reals^d$ with covariance of the form $C_X + R$ where $C_X$ is the covariance of a $\star$-scale invariant field and $R$ is the kernel of an integral operator which maps $L^2(\reals^d)$ to $H^{s}(\reals^d)$ for some $s > d$ (in particular it is in this sense more regular than $C_X$ which maps $L^2(\reals^d)$ to $H^d(\reals^d)$). The second step is then to actually make the calculations work, and to do this in the general set-up we make use of some operator-theoretic methods.

\subsection{Extension of log-correlated fields to the whole space}

Let us begin by solving the aforementioned extension problem.
In what follows we will denote by the same symbols both the integral operators and their kernels, and $C_X$ (resp. $C_{Y^{(\alpha)}}$) will always refer to the covariance operator of a $\star$-scale (resp. almost $\star$-scale) invariant field with a fixed seed covariance $k$ (resp. and parameter $\alpha$).

First of all, we note the existence of the following partition of unity consisting of squares of smooth functions.

\begin{lemma}\label{lemma:partition_of_unity}
  Let $U \subset \reals^d$ be an open domain and $V \Subset U$ an open subdomain.
  Then there exists an open set $W$ with $V \Subset W \Subset U$ and non-negative functions $a,b \in C^\infty(\reals^d)$ such that $a^2 + b^2 \equiv 1$, $b(x) = 0$ for all $x \in \overline{V}$, $b(x) > 0$ for all $x \in \reals^d \setminus \overline{V}$ and $a(x) = 0$ for all $x \in \reals^d \setminus W$. 
\end{lemma}

\begin{proof}
  Pick any $W$ with $V \Subset W \Subset U$.
  It is well-known that one can pick a function $u \in C^\infty(\reals^d)$ which is $1$ in $V$, $0$ outside $W$ and $0 \le u(x) < 1$ for $x \in W \setminus \overline{V}$.
  The function $u(x)^2 + (1 - u(x))^2 \ge \frac{1}{2}$ is everywhere strictly positive and therefore the function $v(x) \coloneqq \sqrt{u(x)^2 + (1 - u(x))^2}$ is smooth and strictly positive.
  Finally define $a(x) \coloneqq u(x)/v(x)$ and $b(x) \coloneqq (1 - u(x))/v(x)$ to obtain the desired functions.
\end{proof}

Secondly we need the following estimates on the covariance operator $C_X$.

\begin{lemma}\label{lemma:C_Y_fourier}
  For any $s \in \reals$ the operator $C_X$ is a bounded invertible operator $H^s(\reals^d) \to H^{s+d}(\reals^d)$.
  The same holds for $C_{Y^{(\alpha)}}$ for any $\alpha > 0$.
  In particular the Cameron--Martin space of $Y^{(\alpha)}$ equals $H^{d/2}(\reals^d)$ with an equivalent norm.
  
  Moreover the Fourier transform of the associated kernel
  \[K(u) \coloneqq C_X(u,0) = \int_0^{\infty} k(e^s u) \, ds\]
  is smooth and satisfies
  \[|\nabla_\xi \hat{K}(\xi)| \lesssim (1 + |\xi|^2)^{-\frac{d+1}{2}}.\]
\end{lemma}

\begin{proof}
  We have $C_X f = K * f$, so it is enough to study the Fourier transform of $K$.
  We compute
  \[\hat{K}(\xi) = \int_0^{\infty} e^{-du} \hat{k}(e^{-u} \xi) du = \int_0^1 v^{d-1} \hat{k}(v \xi) \, dv = |\xi|^{-d} \int_0^{|\xi|}v^{d-1} \hat{k}(v) \, dv.\]
  Since $\hat{k}(0) > 0$ and also $\hat{k}(\xi) = O(|\xi|^{-\alpha})$ for some $\alpha > d+1$, we see that the above quantity is bounded from above and below by a constant multiple of $(1 + |\xi|^2)^{-d/2}$, which implies the claim that $C_X$ maps $H^s(\reals^d)$ to $H^{s+d}(\reals^d)$ continuously and bijectively.
  
  Similarly $C_{Y^{(\alpha)}} f = K_\alpha * f$ with
  \[\hat{K}_{\alpha}(\xi) = \int_0^1 v^{d-1} \hat{k}(v \xi)(1 - v^\alpha) \, dv = |\xi|^{-d} \int_0^{|\xi|} v^{d-1}\hat{k}(v)(1 - |\xi|^{-\alpha} v^\alpha) \, dv\]
  and one again sees that this is bounded from above and below by a constant multiple of $(1 + |\xi|^2)^{-d/2}$.
  In particular $H_{Y^{(\alpha)}} = C_{Y^{(\alpha)}}^{1/2} L^2(\reals^d) = H^{d/2}(\reals^d)$.

  Next we note that since $k$ is compactly supported, $\hat{k}$ is smooth and also $|\nabla \hat{k}(\xi)| = O(|\xi|^{-\alpha})$.
  Thus
  \[\nabla \hat{K}(\xi) = \int_0^1 v^{d} \nabla \hat{k}(v\xi) dv = |\xi|^{-d-1} \int_0^{|\xi|} v^{d} \nabla \hat{k}(v) \, dv,\]
  from which the second claim follows.
\end{proof}

As a corollary of the following lemma from \cite{JSW2} we can rephrase \eqref{eq:logkernel} using a $\star$-scale invariant covariance instead of pure logarithm.

\begin{lemma}[{\cite[Proposition~4.1~(vi)]{JSW2}}]\label{lemma:y_cov_decomposition}
  The covariance $C_X$ of a $\star$-scale invariant field $X$ satisfies $C_X(x,y) = \log \frac{1}{|x-y|} + g_0(x,y)$, where $g_0(x,y)$ belongs to $H^{s'}(\reals^d)$ for some $s' > d$.
\end{lemma}

Let us next prove the extension itself. We emphasise that the kernel $R$ in the proposition below is not necessarily definite positive.

\begin{proposition}\label{lemma:extension}
  Let $C_\Gamma$ be as in Theorem \ref{thm:decomposition}. Let $V \Subset U$ be an open subdomain. Let $X$ be a $\star$-scale invariant log-correlated field with a seed covariance $k$ satisfying Assumption \ref{assumption:seed_covariance}.
  
  Then there exists a bounded integral operator $R \colon L^2(\reals^d) \to L^2(\reals^d)$ such that $C_X + R$ is strictly positive and the corresponding kernels satisfy
  \[C_\Gamma(x,y) = C_X(x,y) + R(x,y)\]
  for all $x,y \in V$.
  The kernel $R$ is Hölder-continuous with some exponent $\gamma > 0$ and moreover, there exists $\delta > 0$ such that $R$ defines a bounded operator $H^{r}(\reals^d) \to H^{r+d+2\delta}(\reals^d)$ for all $r \in [-d,0]$.
\end{proposition}

\begin{proof}
  Let $V \Subset W \Subset U$ and $a,b \in C^\infty(\reals^d)$ be as in Lemma~\ref{lemma:partition_of_unity} and consider the (distribution-valued) Gaussian field $Z = a \Gamma + b X$ defined on $\reals^d$.
  Here $\Gamma$ and $X$ are independent and have covariance operators $C_\Gamma$ and $C_X$ respectively.
  By using Lemma~\ref{lemma:y_cov_decomposition} we can write $C_\Gamma(x,y) = C_X(x,y) + \tilde{g}(x,y)$ with $\tilde{g} \in H^{s'}_{\mathrm{loc}}(\reals^d \times \reals^d)$ for some $s' > d$.
  Thus we may write the kernel of the covariance operator of $Z$ as
  \[C_Z(x,y) = a(x)a(y) C_\Gamma(x,y) + b(x)b(y) C_X(x,y) = C_X(x,y) + R(x,y),\]
  where
  \begin{equation}\label{eq:H}
  R(x,y) \coloneqq (a(x)a(y) + b(x)b(y) - 1)C_X(x,y) + a(x)a(y) \tilde{g}(x,y).
  \end{equation}
  Note that $G(x,y) \coloneqq a(x)a(y)\tilde{g}(x,y)$ is an element of $H^{s'}(\reals^d \times \reals^d)$. For any $f \in H^r(\reals^d)$ with $r \in [-s',0]$ we have that the corresponding operator $G$ satisfies
  \begin{align*}
    \|G f\|_{H^{r + s'}(\reals^d)}^2 & =  \int_{\reals^d} (1 + |\xi|^2)^{r + s'} \Big|\int_{\reals^d} \hat{G}(\xi, \zeta) \overline{\hat{f}(\zeta)} \, d\zeta\Big|^2 \, d\xi \\
    & \lesssim \|G\|_{H^{s'}(\reals^d \times \reals^d)}^2 \|f\|_{H^r(\reals^d)}^2.
  \end{align*}
We conclude that $G$ is a bounded operator $H^{r}(\reals^d) \to H^{r+s'}(\reals^d)$.

  Let us then consider the operator $T$ with kernel
  \[T(x,y) \coloneqq (a(x)a(y) + b(x)b(y) - 1)C_X(x,y)\]
  corresponding to the first term in the definition of $R$.
  Again for $f \in L^2(\reals^d)$ we have
  \[\|T f\|_{H^{d+1}(\reals^d)}^2 = \int_{\reals^d} (1 + |\xi|^2)^{d+1} \Big|\int_{\reals^d} \hat{T}(\xi,\zeta) \overline{\hat{f}(\zeta)} \, d\zeta\Big|^2 \, d\xi.\]
  Note that since $a^2 + b^2 = 1$ we have
  \[T(x,y) = (a(x)(a(y) - a(x)) + b(x)(b(y) - b(x)))C_X(x,y).\]
  The maps $f \mapsto a f$ and $f \mapsto b f = (b - 1)f + f$ are bounded operators $H^{\alpha}(\reals^d) \to H^{\alpha}(\reals^d)$ for any $\alpha \in \reals$ since $a$ and $b-1$ are compactly supported and smooth.
  Thus it is enough to show that $A \colon f \mapsto \big[x \mapsto \int (a(y) - a(x))K(x-y) f(y) \, dy\big]$ and $B \colon f \mapsto \big[x \mapsto \int (b(y) - b(x))K(x-y) f(y) \, dy\big]$ are bounded operators $H^r(\reals^d) \to H^{r+d+1}(\reals^d)$, where $K(u) = C_X(u,0)$.
  
  We will show the claim for $A$ -- the same proof works for $B$ as well since we only use the fact that $a$ is smooth and has compact support and we can again reduce to this situation by replacing $b$ with $b-1$.
  
  The boundedness of $A \colon H^r(\reals^d) \to H^{r+d+1}(\reals^d)$ boils down to showing that for any $f \in H^r(\reals^d)$ we have the inequality
  \begin{equation}\label{eq:Aboundedness}
      \int (1 + |\xi|^2)^{r+d+1} |\widehat{A f}(\xi)|^2 \, d\xi \lesssim \int (1 + |\xi|^2)^{r} |\hat{f}(\xi)|^2 \, d\xi.
  \end{equation}
  
  A small computation shows that we can write
  \begin{equation*}
    \widehat{A f}(\xi) =   \int_{\reals^d} \hat{a}(\xi - \zeta)(\hat{K}(\xi) - \hat{K}(\zeta))\hat{f}(\zeta) \, d\zeta \\
  \end{equation*}
  We can bound 
  \begin{align*}
       \int_{\reals^d} \hat{a}(\xi - \zeta)(\hat{K}(\xi) - \hat{K}(\zeta))\hat{f}(\zeta) \, d\zeta 
    & \lesssim \int_{\reals^d \setminus B(\xi,|\xi|/2)} |\hat{a}(\xi - \zeta)| |\hat{f}(\zeta)| \, d\zeta \\
    & \quad + \int_{B(\xi,|\xi|/2)} |\hat{a}(\xi - \zeta)| |\xi - \zeta| \sup_{z \in B(\xi,|\xi|/2)} |\nabla \hat{K}(z)| |\hat{f}(\zeta)| \, d\zeta.
   \end{align*}
   By using the smoothness of $a$, we have for $\zeta \in \reals^d \setminus B(\xi, |\xi|/2)$ the inequality $|\hat{a}(\xi - \zeta)| \lesssim (1 + |\xi|^2)^{d-1} (1 + |\zeta|^2)^{\frac{r - d - 1}{2}}$.
   By Cauchy--Schwarz we can therefore bound the first term by
   \[\lesssim (1 + |\xi|^2)^{d-1} \Big(\int_{\reals^d} (1 + |\zeta|^2)^{-d-1} \, d\zeta \Big)^{1/2} \Big( \int_{\reals^d} (1 + |\zeta|^2)^{r} |\hat{f}(\zeta)|^2 \, d\zeta \Big)^{1/2} \lesssim (1 + |\xi|^2)^{-d-1} \|f\|_{H^r(\reals^d)}.\]
   This combined with using Lemma~\ref{lemma:C_Y_fourier} to bound the second term we get
   \begin{align*}
     \widehat{A f}(\xi) \lesssim (1 + |\xi|^2)^{-d-1} \|f\|_{H^r(\reals^d)} + (1 + |\xi|^2)^{-\frac{d+1}{2}} \int_{\reals^d} |\hat{a}(\xi - \zeta)| |\xi - \zeta| |\hat{f}(\zeta)| \, d\zeta.
  \end{align*}
  Thus recalling that we want to prove \eqref{eq:Aboundedness} we have
  \begin{align*}
      & \int (1 + |\xi|^2)^{r+d+1} |\widehat{A f}(\xi)|^2 \\ & \lesssim \int (1 + |\xi|^2)^{r+d+1} (1 + |\xi|^2)^{-2d-2} \|f\|_{H^r(\reals^d)}^2 \\
      & + \int (1 + |\xi|^2)^{r+d+1} (1 + |\xi|^2)^{-d-1} \left(\int_{\reals^d} |\hat{a}(\xi - \zeta)||\xi - \zeta| |\hat{f}(\zeta)| d\zeta\right)^2.
  \end{align*}
  Now, as $r < 0$, the first term is bounded by a constant times $\|f\|_{H^r(\reals^d)}^2$.
  For the second term we let $p(\xi) \coloneqq |\xi||\hat{a}(\xi)|$ and note that since $|\hat{f}(\zeta)| |\hat{f}(\zeta')| \le (|\hat{f}(\zeta)|^2 + |\hat{f}(\zeta')|^2)/2$ we have
  \begin{align*}
    & \int_{\reals^d} (1 + |\xi|^2)^{r+d+1}  (1 + |\xi|^2)^{-d-1}  \left(\int_{\reals^d} p(\xi - \zeta) |\hat{f}(\zeta)|  \, d\zeta\right)^2 d\xi \\
    & = \int_{\reals^d} \int_{\reals^d} \int_{\reals^d} (1 + |\xi|^2)^{r} p(\xi - \zeta) p(\xi - \zeta') |\hat{f}(\zeta)| |\hat{f}(\zeta')| \, d\zeta \, d\zeta' \, d\xi \\
    & \le \int_{\reals^d} \int_{\reals^d} \int_{\reals^d} (1 + |\xi|^2)^r p(\xi - \zeta) p(\xi - \zeta') |\hat{f}(\zeta)|^2 \, d\zeta \, d\zeta' \, d\xi.
  \end{align*}
  Integrating over $\zeta'$ gives just $\|p\|_{L^1(\reals^d)}$ and then by using the inequality $(1 + |\xi|^2)^r \lesssim (1+|\zeta - \xi|^2)^{-r} (1+|\zeta|^2)^{r}$ we may also integrate over $\xi$ and $\zeta$ separately to see that the above is bounded by a constant times
  \[\|p\|_{L^1(\reals^d)} \|(1 + |\cdot|)^{-r}p(\cdot)\|_{L^1(\reals^d)} \|f\|_{H^r(\reals^d)}^2.\]
  Thus putting things together we obtain \eqref{eq:Aboundedness}. Overall we have shown that $R$ as defined in \eqref{eq:H} maps $H^{r}(\reals^d) \to H^{r + d + 2\delta}$ for $\delta > 0$ small enough.
 
 Let us next show that $R$ is Hölder-continuous.
 As $\tilde{g}$ belongs to $H^{s'}_{\mathrm{loc}}(\reals^d \times \reals^d)$ for some $s' > d$, it follows from the Sobolev embedding $H^{d + \delta}(\reals^{2d}) \to C^\delta(\reals^{2d})$ where $C^\delta(\reals^{2d})$ is the space of $\delta$-Hölder functions vanishing at infinity, that $\tilde{g}$ is $\gamma$-Hölder for some $\gamma > 0$.
 By \eqref{eq:H} this implies that we only need to show that $(a(x)a(y) + b(x)b(y) - 1)C_X(x,y)$ is Hölder-continuous.
 As this term is compactly supported, we can add a smooth cutoff function $\rho$ such that
 \[(a(x)a(y) + b(x)b(y) - 1)C_X(x,y) = \rho(x)\rho(y)(a(x)(a(y) - a(x)) + b(x)(b(y) - b(x))) C_X(x,y)\]
 for all $x,y \in \reals^d$.
 Moreover, since $C_X(x,y) = \log \frac{1}{|x-y|} + g_0(x,y)$ with $g_0$ smooth, it is enough to show that
 \[(a(y) - a(x))\rho(x)\rho(y)\log \frac{1}{|x-y|}\]
 is Hölder-continuous (the term with $b(y)-b(x)$ can again be handled in a similar manner).
 Let us write the above as
 \[\int_0^1 \nabla a(x + u(y-x)) \, du \cdot (y-x) \rho(x)\rho(y) \log \frac{1}{|x-y|}.\]
 As $a$ is smooth, the map $(x,y) \mapsto \int_0^1 \nabla a(x + u(y-x)) \, du$ is in particular a Hölder continuous map $\reals^{2d} \to \reals^d$.
 Thus it is enough to show that
 $(x,y) \mapsto (y-x) \log \frac{1}{|x-y|}$
 is Hölder-continuous but this follows easily by checking that each component function $(y_j - x_j) \log \frac{1}{|x-y|}$ is Hölder continuous in each coordinate.
 The Hölder constants are also easily seen to be bounded for $x,y \in \supp \rho$.

  Finally let us note that $C_Z$ is strictly positive since if $f \in L^2(\reals^d)$ is nonzero, then at least one of $f|_{V}$ or $f|_{\supp b}$ is nonzero.
  In the first case $\int a(x)a(y)C_\Gamma(x,y) f(x)f(y) > 0$ by the assumption that $C_\Gamma$ was assumed to be injective in $V$, while in the second case $\int b(x)b(y)C_X(x,y) f(x)f(y) > 0$ since $C_X$ is strictly positive on whole of $\reals^d$.
\end{proof}

\subsection{Deducing the decomposition theorem}

Having obtained the desired extension, we are ready to prove the decomposition theorem. The second part of the proof consists in showing that we may subtract $C_{Y^{(\alpha)}}$ from $C_X + R$ for some small enough $\alpha > 0$ and still obtain a positive operator.

To do this, we need to use the following classical stability property of strictly positive operators of the form $1 + K$ with $K$ compact and self-adjoint that follows directly from the spectral theorem.

\begin{lemma}\label{lemma:specstab}
Let $\mathcal{H}$ be a Hilbert space and $T$ a self-adjoint compact operator on $\mathcal{H}$ and suppose that $1 + T$ is strictly positive. Then there exists $\varepsilon > 0$ such that $1 + A + T$ is strictly positive for any self-adjoint $A$ with $\|A\|_{\mathcal{H} \to \mathcal{H}} \leq \varepsilon$. 
\end{lemma}

As a consequence of the above lemma and the smoothing properties of the map $R$ obtained in Lemma~\ref{lemma:extension} we first create a necessary lee-room. Notice that $C_X + R = C_X^{1/2}(I + C_X^{-1/2}RC_X^{-1/2})C_X^{1/2}$ and hence
 $$\langle (C_X + R)f, f \rangle_{L^2(\R^d)} = \langle (I + C_X^{-1/2}RC_X^{-1/2})C_X^{1/2}f, C_X^{1/2}f \rangle_{L^2(\R^d)}.$$
 The following statement is thus effectively saying that in fact $\langle (C_X + R)f, f \rangle_{L^2(\R^d)} > 0$ not only for $f \in L^2(\R^d)$, but also for $f \in H^{-d/2}(\reals^d)$.
 
\begin{lemma}\label{prop:epspert}
There is some $\varepsilon > 0$ such that $1 + A + C_X^{-1/2} R C_X^{-1/2}$ is a strictly positive operator on $L^2(\reals^d)$ for any self-adjoint $A$ with $\|A\|_{L^2(\reals^d) \to L^2(\reals^d)} \le \varepsilon$.
\end{lemma}

\begin{proof}
  We start by observing that the operator $\tilde{R} = C_X^{-1/2} R C_X^{-1/2}$ is compact from $L^2(\R^d)$ to $L^2(\R^d)$. Indeed, we can write $\tilde{R}$ as $C_X^{-1/2} J R C_X^{-1/2}$ where $J$ is the identity map. Now, due to the fact that $R(x,y)$ has compact support (see Equation \eqref{eq:H} and recall that $C_X(x,y) = 0$ for $|x-y| > 1$) this mapping takes successively $$L^2(\R^d) \rightarrow H^{-d/2}(\R^d) \rightarrow H^{d/2 + 2\delta}(B) \rightarrow H^{d/2}(B) \rightarrow L^2(\R^d),$$ where $B \subset \reals^d$ is some fixed large enough open ball such that $B \times B \supset \supp R$. The identity map $J$ from $H^{d/2 + 2\delta}(B) \rightarrow H^{d/2}(B)$ is compact by Rellich-Kondrachov theorems for fractional Sobolev spaces (see e.g. Chapters 1, 2 in \cite{Triebel}) and as the other maps are bounded, the whole composition is compact.
  
  As $R$ is also self-adjoint on $L^2(\reals^d)$, there is an orthonormal basis of $L^2(\reals^d)$ consisting of eigenfunctions of $\tilde{R}$.
  To show that $1 + \tilde{R}$ is strictly positive it is enough to show that $\tilde{R}$ has no eigenfunctions with eigenvalues $\le -1$.
  Assume that $f$ is an eigenfunction of $\tilde{R}$ with nonzero eigenvalue $\lambda$.
  Then by Lemma~\ref{lemma:extension} we know that $\tilde{R}$ maps $H^{s}(\reals^d) \to H^{s + 2\delta}(\reals^d)$ for any $s \in [0,d/2]$ and thus after applying $\tilde{R}$ to $f$ roughly $1/\delta$ times we see that actually $f \in H^{d/2}(\reals^d)$.
  Thus there exists some $g \in L^2(\reals^d)$ such that $f = C_X^{1/2} g$, and we have that
  \[(1 + \lambda) \|f\|_{L^2(\reals^d)}^2 = \langle (1 + \tilde{R})f, f\rangle_{L^2(\reals^d)} = \langle (1 + \tilde{R}) C_X^{1/2} g, C_X^{1/2} g \rangle_{L^2(\reals^d)} = \langle (C_X + R) g, g \rangle_{L^2(\reals^d)} > 0\]
  by the assumption on $C_X+R$, implying that $\lambda > -1$.
  Thus $1 + \tilde{R}$ is strictly positive and the claim follows from Lemma~\ref{lemma:specstab}.
\end{proof}

The final important technical ingredient is that for any $\alpha_0 > 0$, 
$$(C_X - C_{Y^{(\alpha)}})^{-1/2} - C_X^{-1/2} \colon L^2(\reals^d) \to H^{\frac{-d-\alpha_0}{2}}(\reals^d)$$ converges pointwise to $0$ when we let the parameter $\alpha$ of the almost $\star$-scale invariant field $Y^{(\alpha)}$ to $0$.

\begin{lemma}\label{lemma:fourier_estimates}
	For all $\alpha > 0$ set $U_\alpha \coloneqq C_X - C_{Y^{(\alpha)}}$ and let $U_0 = C_X$.
  Then $U_\alpha^{1/2}$ is a bounded bijection $H^{s}(\reals^d) \to H^{s + \frac{d + \alpha}{2}}(\reals^d)$ for all $s \in \reals$, and for any $\alpha_0 > 0$, we have $$\sup_{\alpha_0 \geq \alpha > 0} \|U_{\alpha}^{-1/2}\|_{L^2(\reals^d) \to H^{-\frac{d+\alpha_0}{2}}(\reals^d)} < \infty.$$
  Moreover, for any fixed $\alpha_0>0$ and $f \in L^2(\reals^d)$ we have
  \[\lim_{\alpha \to 0} \|(U_{\alpha}^{-1/2} - C_{Y^{(\alpha)}}^{-1/2}) f\|_{H^{-\frac{d+\alpha_0}{2}}(\reals^d)} = 0.\]
\end{lemma}

Before proving the lemma, let us see how it implies the theorem:

\begin{proof}[Proof of Theorem \ref{thm:decomposition}:]
	We begin by writing
	\[\langle (C_X - C_{Y^{(\alpha)}} + R)f, f\rangle_{L^2(\reals^d)} = \langle (1 + \tilde{R}_\alpha)U_\alpha^{1/2} f, U_\alpha^{1/2} f\rangle_{L^2(\reals^d)},\]
	where $U_\alpha = C_X - C_{Y^{(\alpha)}}$ and $\tilde{R}_\alpha = U_\alpha^{-1/2} R U_\alpha^{-1/2}$. It thus suffices to show that for some sufficiently small $\alpha > 0$ we have
	\[\langle (1 + \tilde{R}_\alpha) g, g \rangle_{L^2(\reals^d)} > 0\]
	for all nonzero $g \in L^2(\reals^d)$.
	Indeed, this implies that $C_X - C_{Y^{(\alpha)}} + R$
	is a positive integral operator on $L^2(\reals^d)$, whose kernel by Lemma~\ref{lemma:extension} and \cite[Proposition~4.1 (iii)]{JSW2} is Hölder-continuous, and thus the corresponding Gaussian process has an almost surely Hölder-continuous version (see e.g. \cite[Theorem~1.3.5]{AdlerTaylor}).
	In addition by Lemma~\ref{lemma:extension} and Lemma~\ref{lemma:fourier_estimates} we see that $R$ and $C_X - C_{Y^{\alpha}}$ map $H^s(\reals^d) \to H^{s+d+\varepsilon}(\reals^d)$ for some $\varepsilon > 0$ and all $s \in [-d,0]$.
	
	To show that $1 + \tilde{R}_\alpha$ is positive on $L^2(\reals^d)$ on the other hand we may write $1 + \tilde{R}_\alpha = 1 + \tilde{R} + (\tilde{R}_\alpha - \tilde{R})$, where $\tilde{R} = C_X^{-1/2} R C_X^{-1/2}$.
	By Lemma~\ref{prop:epspert} it is enough to show that $\|\tilde{R}_\alpha - \tilde{R}\|_{L^2(\reals^d) \to L^2(\reals^d)}$ can be made as small as we wish by choosing $\alpha$ small.
	
	As $\tilde{R}_\alpha - \tilde{R}$ is self-adjoint we have
	$$\|\tilde{R}_\alpha - \tilde{R}\|_{L^2(\reals^d) \to L^2(\reals^d)} = \sup_{u \in L^2(\reals^d), ||u||_2 = 1} |\langle (\tilde{R}_\alpha - \tilde{R}) u, u \rangle|_{L^2(\reals^d)}.$$
	By linearity and self-adjointness of $C_X^{-1/2}, R$ and $U_\alpha^{-1/2}$, we can write $ \langle (\tilde{R}_\alpha - \tilde{R} )u, u \rangle_{L^2(\reals^d)}$ as 
	$$\langle  (U_\alpha^{-1/2} - C_X^{-1/2})R C_X^{-1/2}u, u \rangle_{L^2(\reals^d)} + \langle (U_\alpha^{-1/2} - C_X^{-1/2})R U_\alpha^{-1/2}u, u \rangle_{L^2(\reals^d)}.$$
	Now choose $\alpha_0 = \delta$ in Lemma \ref{lemma:fourier_estimates} and observe that then for all $\alpha < \alpha_0$, the unit ball of $L^2(\reals^d)$ under $R U_\alpha^{-1/2}$ and $R C_X^{-1/2}$ is contained in a fixed compact set of $H^{\frac{d+\delta}{2}}(\reals^d)$. As Lemma \ref{lemma:fourier_estimates} establishes uniform boundedness as well as pointwise convergence, we have that $U_\alpha^{-1/2} \to C_X^{-1/2}$ uniformly on this set and thus conclude the theorem.
\end{proof}

We finally prove the lemma:

\begin{proof}[Proof of Lemma \ref{lemma:fourier_estimates}]
	Note that $U_\alpha$ is a Fourier multiplier operator with the symbol
	\[\hat{u}_{\alpha}(\xi) = \int_0^1 v^{d-1+\alpha} \hat{k}(v\xi) \, dv =  |\xi|^{-d-\alpha} \int_0^{|\xi|} v^{d-1+\alpha} \hat{k}(v) \, dv.\]
	As by assumption $\hat k$ is non-negative and decays faster than any polynomial, we have that 
	\[(1 + |\xi|^2)^{-\frac{d+\alpha}{2}} \lesssim \hat{u}_{\alpha}(\xi) \lesssim (1 + |\xi|^2)^{-\frac{d+\alpha}{2}}\]
	where the hidden constant does not depend on $\alpha$. In particular for every $\alpha < \alpha_0$, we have $(1 + |\xi|^2)^{-\frac{d+\alpha_0}{2}} \lesssim \hat{u}_{\alpha}(\xi)$.

  Let us now fix $\alpha_0$ and consider for $\alpha < \alpha_0$ the self-adjoint operator $T_\alpha = U_\alpha^{-1/2} - C_Y^{-1/2}$ which maps $L^2(\reals^d)$ to $H^{-\frac{d+\alpha}{2}}(\reals^d) \subseteq H^{-\frac{d+\alpha_0}{2}}(\reals^d)$. 
  For any fixed $f \in L^2(\reals^d)$ we have
  \[\|T_\alpha f\|_{H^{-\frac{d+\alpha_0}{2}}(\reals^d)} = \int_{\reals^d} (1+|\xi|^2)^{-\frac{d+\alpha_0}{2}}|\hat{u}_\alpha(\xi)^{-1/2} - \hat{K}(\xi)^{-1/2}|^2 |\hat{f}(\xi)|^2 \, d\xi.\]
  For any fixed $\xi$ the integrand tends to $0$ as $\alpha \to 0$. Thus, as $\hat{u}_\alpha(\xi) \gtrsim (1 + |\xi|^2)^{-\frac{d + \alpha_0}{2}}$ for all $\alpha< \alpha_0$, we can apply the dominated convergence theorem to deduce that $T_\alpha f \to 0$ in $H^{-\frac{d+\alpha_0}{2}}(\reals^d)$.
\end{proof}

\section{General bounds on $\det_\gamma M$ and $\delta(A)$}\label{sec:mvclc}

In this section we prove two (to our knowledge) non-standard lemmas for Malliavin calculus, that we believe could possibly be of independent interest for proving the existence of density and its positivity also in more general settings. Firstly, we prove a certain projection bound for the determinant of complex Malliavin variables. Second, we obtain an estimate on the complex covering fields that is again a much easier starting point for further calculations.

\subsection{Proof of the projection bound -- Proposition \ref{lem:deltalowerbound}}

\begin{proof}[Proof of Proposition \ref{lem:deltalowerbound}]
	Let us first expand
	\begin{align*}
	& \|DF\|^2_{\HC}\Big\|DF - \frac{\langle DF, D\overline{F}\rangle_{\HC}}{\|DF\|_{\HC}^2} D\overline{F}\Big\|_{\HC}^2 \\
	& = \|DF\|^2_{\HC} \Big( \|DF\|^2_{\HC} - \frac{\overline{\scalar{DF,D\overline{F}}}_{\HC}}{\|DF\|^2_{\HC}} \scalar{DF,D\overline{F}}_{\HC} \\
	& ~~~~~- \frac{\scalar{DF,D\overline{F}}_{\HC}}{\|DF\|^2_{\HC}} \scalar{D \overline{F},DF}_{\HC} + \frac{|\scalar{DF,D\overline{F}}_{\HC}|^2}{\|DF\|^4_{\HC}} \|D\overline{F}\|^2_{\HC} \Big) \\
	& = \|DF\|^4_{\HC} - |\scalar{DF,D\overline{F}}_{\HC}|^2.
	\end{align*}
	By \eqref{eq:malliavin_det}, we deduce that
	\begin{equation}\label{E:lastE}
\det \gamma_F = \frac{1}{4} \|DF\|^2_{\HC}\Big\|DF - \frac{\langle DF, D\overline{F}\rangle_{\HC}}{\|DF\|_{\HC}^2} D\overline{F}\Big\|_{\HC}^2.
	\end{equation}
	As we have the following projection inequality
	\[\|DF\|_{\HC} \ge \Big\|DF - \frac{\langle DF, D\overline{F}\rangle_{\HC}}{\|DF\|_{\HC}^2} D\overline{F}\Big\|_{\HC},\]
	the result follows, once we show that for any $h \in \HC$,
	\begin{equation}\label{eq:proj}
	\Big\|DF - \frac{\langle DF, D\overline{F}\rangle_{\HC}}{\|DF\|_{\HC}^2} D\overline{F}\Big\|_{\HC} \ge \frac{\big||\langle DF, h \rangle_{\HC}| - |\langle D\overline{F}, h \rangle_{\HC}|\big|}{\|h\|_{\HC}}.    
	\end{equation}
	By Cauchy--Schwarz inequality and the triangle inequality we have
	\begin{align*}
	\Big\|DF - \frac{\langle DF, D\overline{F}\rangle_{\HC}}{\|DF\|_{\HC}^2} D\overline{F}\Big\|_{\HC} & \ge \frac{|\langle DF - \frac{\langle DF, D\overline{F}\rangle_{\HC}}{\|DF\|_{\HC}^2} D\overline{F}, h\rangle_{\HC}|}{\|h\|_{\HC}} \\
	& \ge \frac{|\langle DF, h \rangle_{\HC}| - \frac{|\langle DF, D\overline{F}\rangle_{\HC}|}{\|DF\|_{\HC}^2} |\langle D\overline{F}, h \rangle_{\HC}|}{\|h\|_{\HC}} \\
	& \ge \frac{|\langle DF, h\rangle_{\HC}| - |\langle D\overline{F}, h\rangle_{\HC}|}{\|h\|_{\HC}}.
	\end{align*}
	By now repeating the bound with $\overline{h}$ in place of $h$ we obtain \eqref{eq:proj} which finishes the proof.
\end{proof}

\subsection{Bounding $\delta(A)$ via derivatives in independent Gaussian directions -- Proposition \ref{prop:bndcv}}

For a succinct write-up, it is helpful to use directional derivatives in independent random directions, although the proposition could also be proved by first proving a version for smooth random variables and then taking limits.

Now, recall that for smooth random variables $F$, and $h \in H_\C$  we could write
\begin{equation}\label{eq:frechet}
  \langle D F(\Gamma), h\rangle_H = \frac{d}{dt}\Big|_{t=0} F(\Gamma + t h).
\end{equation}

We consider directional derivatives in independent random directions, with the law of $\Gamma$. More precisely, let $X \sim \Gamma$ be an independent Gaussian field defined on a new probability space $(\Omega_X,\mathcal{F}_X,\P_X)$ whose expectation we denote by $\E_X$. For a Malliavin variable $F \in \D^{2,\infty}$, as $DF \in H_\C$ and $X$ is independent of $\Gamma$, one can define 
\begin{equation}\label{eq:randdir}
\Dc_X F \coloneqq \langle X, D F(\Gamma) \rangle_H
\end{equation}
and directly conclude from this definition that:

\begin{lemma}\label{lemma:averaging}
Let $X \sim \Gamma$ be independent of $\Gamma$ and $F,G \in \D^{1,\infty}$. We then have that $\E_X[\Dc_X F \cdot \overline{\Dc_X G}] = \langle D F, D G \rangle_{\HC}$.
\end{lemma}

We are now ready to prove Proposition \ref{prop:bndcv}.

\begin{proof}[Proof of Proposition \ref{prop:bndcv}]
Write $\Delta := 4\det \gamma_F = \|DF\|_{\HC}^4 - |\langle DF, D\overline{F}\rangle_{\HC}|^2$. Then by the integration by parts rule for the divergence operator $\delta$ (e.g. \cite[Proposition 1.3.3]{Nualart}), $\delta(A)$ equals
\[\frac{\|DF\|_{\HC}^2 \delta(DF) - \langle DF, D\overline{F}\rangle_{\HC} \delta(D\overline{F})}{\Delta} - \langle D\frac{\|DF\|_{H_{\complexes}}^2}{\Delta}, D\overline{F}\rangle_{H_{\complexes}} + \langle D\frac{\langle DF, D\overline{F}\rangle_{\HC}}{\Delta}, DF\rangle_{\HC}.\]
The first term is $\lesssim \Delta^{-1} \|DF\|_{\HC}^2 |\delta(DF)|$ in absolute value, so it is enough to consider the other two terms. By the product rule for Malliavin derivatives, we may write
\[ \langle D\frac{\langle DF, D\overline{F}\rangle_{\HC}}{\Delta}, DF\rangle_{\HC} - \langle D\frac{\|DF\|_{H_{\complexes}}^2}{\Delta}, D\overline{F}\rangle_{H_{\complexes}} \]
as
\begin{align*}
  & = \Delta^{-1}\left(\langle D\langle DF, D\overline{F}\rangle_{\HC}, DF\rangle_{\HC} - \langle D\|DF\|_{H_{\complexes}}^2, D\overline{F}\rangle_{H_{\complexes}}\right) - \\
  &- \Delta^{-2}\left(\langle DF, D\overline{F}\rangle_{\HC}\langle D \Delta,  DF\rangle_{\HC} - \|DF\|^2_{\HC}\langle  D\Delta, D\overline{F}\rangle_{\HC}\right)
\end{align*}
To bound the first term, we first notice that by Cauchy--Schwarz
\[ \langle D\langle DF, D\overline{F}\rangle_{\HC}, DF\rangle_{\HC} \leq \|D\langle DF, D\overline{F}\rangle_{\HC}\| \|DF\|_{\HC}.\]
For the first term, it is now helpful to use the averaging in Lemma~\ref{lemma:averaging} for a quick bound. We write
\[ \|D\langle DF, D\overline{F}\rangle_{\HC}\|_{\HC} = 2|\E_{X,Y} \Dc_Y F \cdot \Dc_X \Dc_Y F|.\]
By Cauchy--Schwarz this can be bounded by
\[ 2\sqrt{\E_{X,Y} |\Dc_Y F|^2}\sqrt{\E_{X,Y}|\Dc_X \Dc_Y F|^2} = 2\|DF\|_{\HC}\|D^2F\|_{\HC\otimes \HC}.\]
Similarly, one can bound
\[ \langle D\|DF\|_{H_{\complexes}}^2, D\overline{F}\rangle_{H_{\complexes}} \leq 2\|DF\|_{\HC}\|D^2F\|_{\HC\otimes \HC},\]
and thus
\[\Delta^{-1}\left(\langle D\langle DF, D\overline{F}\rangle_{\HC}, DF\rangle_{\HC} - \langle D\|DF\|_{H_{\complexes}}^2, D\overline{F}\rangle_{H_{\complexes}}\right) \leq 4\frac{\|DF\|^2_{\HC}\|D^2F\|_{\HC\otimes \HC}}{\Delta}.\]
It remains to handle 
\[\Delta^{-2} \left( \langle DF, D\overline{F}\rangle_{\HC}\langle D \Delta,  DF\rangle_{\HC} - \|DF\|^2_{\HC}\langle  D\Delta, D\overline{F}\rangle_{\HC} \right),\]
which we can rewrite as
\[\Delta^{-2}\langle D \Delta,  \langle D\overline{F}, DF\rangle_{\HC} DF -\|DF\|^2_{\HC} D\overline{F}\rangle_{\HC}.\]
By Cauchy--Schwarz this expression is bounded by
\[\Delta^{-2} \|D\Delta\|_{\HC} \|\langle D\overline{F}, DF\rangle_{\HC} DF - \|DF\|_{\HC}^2 D\overline{F}\|_{\HC} = \Delta^{-3/2} \|D\Delta\|_{\HC} \|DF\|_{\HC},\]
where we have used the fact (derived in Equation \eqref{E:lastE}) that
\begin{equation}\label{eq:det}
\|DF\|_{\HC}^2 \Delta = \|\langle D\overline{F},DF\rangle_{\HC} DF - \|DF\|^2_{\HC} D\overline{F}\|_{\HC}^2.
\end{equation}
Thus the proposition follows from the following claim:
\begin{claim}
We have that $\|D\Delta\|_{\HC} \lesssim \Delta^{1/2}\|DF\|_{\HC}\|D^2F\|_{\HC\otimes \HC}$. 
\end{claim}

\begin{proof}[Proof of claim]
Maybe the nicest way to prove this claim is to use derivatives in random directions as above. First, observe that using averaging we can write a neat analogue of Equation \eqref{eq:det} :
\[\Delta = \frac{1}{2} \E_{Z,W} |\Dc_Z F \cdot \Dc_W \overline{F} - \Dc_Z \overline{F} \cdot \Dc_W F|^2.\]
Thus we have
\[\Dc_X \Delta = \Re \E_{Z,W} (\Dc_Z F \cdot \Dc_W \overline{F} - \Dc_Z \overline{F} \cdot \Dc_W F)\Dc_X(\Dc_Z F \cdot \Dc_W \overline{F} - \Dc_Z \overline{F} \cdot \Dc_W F).\]
By triangle inequality and Cauchy--Schwarz we obtain
\[|\Dc_X \Delta|^2 \lesssim \Delta \E_{Z,W} |\Dc_X(\Dc_Z F \cdot \Dc_W \overline{F})|^2\]
and hence
\[\|D\Delta\|_{\HC}^2 = \E_X |\Dc_X \Delta|^2 \lesssim \Delta \|DF\|_{\HC}^2 \|D^2 F\|_{\HC \otimes \HC}^2,\]
from which the claim follows.
\end{proof}
\end{proof}

\section{Estimates for Malliavin variables in the case of imaginary chaos}\label{sec:estimates_chaos}

The aim of this section is to prove the probabilistic bounds needed to apply the tools of Malliavin calculus to $M = \mu(f)$. We start by going through some old and new Onsager inequalities and related integral bounds. In Section \ref{sec:lemma:dinfty}, we prove by a rather standard argument that $M$ is in $\D^\infty$, i.e. Proposition \ref{prop:dinfty}. In Section \ref{sec:bounds_for_derivatives} we derive bounds on $|\delta(DM)|$ and 
$\|D^2 M\|_{H_\C \otimes H_\C}$ and deduce Proposition \ref{prop:bound_delta_D2_intro}
by a quite similar argument.

Finally, in Section \ref{sec:mainproof} we prove bounds on the Malliavin determinant of $M$ and this is the main technical input of the paper. Here things get quite interesting -- we rely both on the decomposition theorem, Theorem \ref{thm:decomposition}, and projection bounds for Mallivan determinants from Section \ref{sec:mvclc}, but also need to find ways to get a good grip on the concentration of $M = \mu(f)$, and on Sobolev norms of the imaginary chaos $\mu$ itself. 

\subsection{Onsager inequalities and related bounds}\label{sec:onsager}

In this section, we collect a few Onsager inequalities and related bounds. To this end, we define for any Gaussian field $\Gamma$ and $\mathbf{x} = (x_1,\dots,x_N), \mathbf{y} = (y_1,\dots,y_M)$ the quantity
\[
\mathcal{E}(\Gamma; \mathbf{x}; \mathbf{y}) = -\sum_{1 \le j < k \le N} \E \Gamma(x_j)\Gamma(x_k) - \sum_{1 \le j < k \le M} \E \Gamma(y_j)\Gamma(y_k) + \sum_{\substack{1 \le j \le N \\ 1 \le k \le M}} \E \Gamma(x_j)\Gamma(y_k).
\]
Also, we let $\Gamma_\delta = \Gamma * \varphi_\delta$ be a mollification of $\Gamma$ where $\varphi_\delta = \delta^{-d} \varphi(\cdot/\delta)$ and $\varphi$ is a smooth non-negative function with compact support that satisfies $\int_{\R^d} \varphi = 1$.

The following is a restatement of a standard Onsager inequality from \cite{JSW}.\footnote{In fact, the cited result does not contain the case of the circle, however essentially the same proof works.}

\begin{lemma}[Proposition 3.6(ii) of \cite{JSW}]\label{lem:Onsager_weak}
Let $K$ be a compact subset of $U$ or the circle $K = S^1$. There exists $C = C(K) >0$ such that the following holds true: Let $N \geq 1, \delta >0$ and for all $i=1 \dots N$ let $x_i, y_i \in K$ be such that $D(x_i, \delta)$ and $D(y_i,\delta)$ are included in $K$. For all $i=1 \dots N$, denote $z_i \coloneqq x_i$ and $z_{N+i} \coloneqq y_i$ and set $d_j \coloneqq \min_{k \neq j} |z_k - z_j|$. Then
\begin{equation}
\label{eq:lem_Onsager_conv}
\Ec(\Gamma_\delta ; \mathbf{x} ; \mathbf{y}) \leq \frac{1}{2} \sum_{j=1}^{2N} \log \frac{1}{d_j} + C N^2.
\end{equation}
Moreover, the same holds for the field $\Gamma$ itself.
\end{lemma}

We will also need stronger Onsager inequalities for (almost) $\star$-scale invariant fields, whose rather standard proof is pushed to the appendix \ref{appendixA}.

\begin{lemma}\label{lemma:star_scale_onsager}
Let $Y_\varepsilon$ and $\hat{Y}_\varepsilon$ be defined as in Section~\ref{sec:starscale} and let $\mathbf{x} = (x_1,\dots,x_N)$ and $\mathbf{y} = (y_1,\dots,y_N)$ be two $N$-tuples of points in $U$.
For all $j = 1,\dots,N$, denote $z_j \coloneqq x_j$ and $z_{N+j} = y_j$ and set $d_j \coloneqq \min_{k \neq j} |z_k - z_j|$.
Then
\[\mathcal{E}(Y_\varepsilon;\mathbf{x};\mathbf{y}) \le \frac{1}{2} \sum_{j=1}^{2N} \log \frac{1}{d_j \vee \varepsilon}\]
and
\begin{equation}
    \label{eq:onsager_tail}
\mathcal{E}(\hat{Y}_\varepsilon(\varepsilon \cdot);\mathbf{x};\mathbf{y}) \le \frac{1}{2} \sum_{j=1}^{2N} \log \frac{1}{d_j}.
\end{equation}
Moreover, if $R$ is a Gaussian field such that $M \coloneqq \sup_{x \in U} \E[R(x)^2] < \infty$, then
\begin{equation}
    \label{eq:onsager_smooth}
\mathcal{E}(R;\mathbf{x};\mathbf{y}) \le N M.
\end{equation}
\end{lemma}

Both of these Onsager inequalities are used in conjunction with the following bounds:

\begin{lemma}\label{lem:integral_min} For $N \geq 2$, there exists $C > 0$ such that

\begin{itemize}
    \item for all $\beta \in (0,\sqrt{d})$,
\begin{equation}\label{eq:lem_bound_integral_1}
    \int_{B(0,1)^N} \prod_{i=1}^N \left( \min_{j \neq i} |z_i - z_j| \right)^{-\beta^2/2} dz_1 \dots dz_N \leq C^N(d-\beta^2)^{-\floor{N/2}} N^{\frac{N\beta^2}{2d}};
\end{equation}

\item for all $\beta \in (0,\sqrt{d})$,
\begin{equation}
\label{eq:lem_int1}
\int_{B(0,1)^N} \prod_{i=1}^N \abs{\log \min_{j \neq i} |z_i - z_j|}^{1/2} \left( \min_{j \neq i} |z_i - z_j| \right)^{-\beta^2/2} d z_1 \dots d z_N
\leq C^N (d-\beta^2)^{- 2 \floor{N/2} } N^N;
\end{equation}

\item for all $\beta \in (0,\sqrt{d})$,
\begin{equation}
\label{eq:lem_int2}
\int_{B(0,1)^N} \prod_{i=1}^N \abs{\log \min_{j \neq i} |z_i - z_j|} \left( \min_{j \neq i} |z_i - z_j| \right)^{-\beta^2/2} d z_1 \dots d z_N
\leq C^N (d-\beta^2)^{- 3 \floor{N/2} } N^N;
\end{equation}

\item for all $\beta > 0$,
\begin{equation}\label{eq:super_critical_blowup}
\int_{B(0,1)^{N}} \left(\prod_{i=1}^{N} \min_{j \neq i} \max(\delta, |z_i - z_j|) \right)^{-\beta^2/2} \, dz_1 \dots dz_N \le C^N N^{N} (\log \frac{1}{\delta})^{N/2} \delta^{-\max(0,\beta^2 - d)N/2};
\end{equation}
\end{itemize}

\end{lemma}

\begin{proof}
We only sketch the proof, as all the main ideas can be found in proof of \cite[Lemma 3.10]{JSW}.

Let us start with showing \eqref{eq:lem_bound_integral_1}. By carefully following the proof of \cite[Lemma~3.10]{JSW} which shows that \eqref{eq:lem_bound_integral_1} is less than $c^{2 \floor{N/2}} N^{\frac{N \beta^2}{2d}}$, one can actually see that the constant $c$ there can be taken to be equal to $c' (d-\beta^2)^{-1/2}$ for some constant $c'>0$ independent of $\beta$ (at one point in the proof there is a term of order $(d-\beta^2)^{-k}$ coming from $\Gamma(1 - \frac{d}{\beta^2})^k$ where $k \le \lfloor N / 2 \rfloor$).

We will next show \eqref{eq:super_critical_blowup}.
By mimicking the beginning of the proof of \cite[Lemma~3.10]{JSW}, we can bound the left hand side of \eqref{eq:super_critical_blowup} by
\[
C^N
\sum_{k=1}^{\lfloor N/2 \rfloor} \sum_F \int_{B(0,1)^N} \prod_{i=1}^k (\delta \vee |u_{2i-1}|)^{-\beta^2} \prod_{i=2k+1}^{N} (\delta \vee |u_i|)^{-\beta^2/2} du_1 \dots du_N
\]
where $C>0$ and the second sum runs over all nearest neighbour configurations $F$ such that the induced graph with vertices $\{1,\dots,N\}$ and edges $(i,F(i))$ has $k$ components. Of course, the domain on which we integrate is actually much smaller than $B(0,1)$, but integrating over this larger domain will be enough for our purposes. After integration, we obtain that the left hand side of \eqref{eq:super_critical_blowup} is at most
\begin{align*}
& C^N \sum_{k=1}^{\lfloor N/2 \rfloor} \sum_F A_{\beta^2}^k A_{\beta^2/2}^{N-2k} \le C^N N^{N} \sum_{k=1}^{\lfloor N/2 \rfloor} A_{\beta^2}^k A_{\beta^2/2}^{N-2k},
\end{align*}
where
\[A_{\beta^2} \coloneqq \int_0^1 r^{d-1} (\delta \vee r)^{-\beta^2} \, dr.\]
Now, by Jensen's inequality $A_{\beta^2/2}^2 \le d^{-1} A_{\beta^2}$, giving us the bound
$C^N N^{N} A_{\beta^2}^{N/2}$. Noting that
\[A_{\beta^2} \lesssim \log \frac{1}{\delta} \delta^{-\max(0,\beta^2 - d)}\]
concludes the proof of \eqref{eq:super_critical_blowup}.

We finally turn to the proof of \eqref{eq:lem_int1} and \eqref{eq:lem_int2}. By again mimicking the beginning of the proof of \cite[Lemma 3.10]{JSW}, we can bound the left hand side of \eqref{eq:lem_int1} by
\begin{align*}
& C^N \sum_{k=1}^{\floor{N/2}} M_k \int_{B(0,1)^N} \prod_{i=1}^k |u_{2i-1}|^{-\beta^2} |\log |u_{2i-1}|| \prod_{i=2k+1}^N |u_i|^{-\beta^2/2} |\log |u_i||^{1/2} \\
& \leq C^N \sum_{k=1}^{\floor{N/2}} M_k \left( \int_0^1 r^{-\beta^2 + d -1} |\log r| d r \right)^k
\leq C^N \sum_{k=1}^{\floor{N/2}} M_k (d-\beta^2)^{-2k} \leq C^N (d-\beta^2)^{-2 \floor{N/2}} N^N,
\end{align*}
where $M_k$ is the number of nearest neighbour functions $\{1,\dots,N\} \to \{1,\dots,N\}$ with $k$ components and $C$ is some large enough constant.
This concludes the proof of \eqref{eq:lem_int1}; the proof of \eqref{eq:lem_int2} is similar.
\end{proof}

\subsection{$M$ belongs to \texorpdfstring{$\D^\infty$}{D infinity} -- proof of Proposition~\ref{prop:dinfty}}\label{sec:lemma:dinfty}
The purpose of this section is to prove Proposition~\ref{prop:dinfty}. Before doing so, we collect two auxiliary lemmas from Malliavin calculus.

\begin{lemma}[{\cite[Lemma 1.2.3]{Nualart}}]\label{lem:Nualart0}
Let $(F_n,n \geq 1)$ be a sequence of (complex) random variables in $\D^{1,2}$ that converges to $F$ in $L^2(\Omega)$ and such that $\sup_n \Expect{ \norme{DF_n}_{H_\C}^2 } < \infty$. Then $F$ belongs to $\D^{1,2}$ and the sequence of derivatives $(DF_n, n \geq 1)$ converges to $DF$ in the weak topology of $L^2(\Omega;H_\C)$.
\end{lemma}

Second, we need a rather direct consequence of \cite[Lemma 1.5.3]{Nualart}:

\begin{lemma}\label{lem:Nualart}
Let $p > 1$, $k \geq 1$ and let $(F_n,n \geq 1)$ be a sequence of (complex) random variables converging to $F$ in $L^p(\Omega)$. Suppose that $\sup_n \norme{F_n}_{k,p} < \infty$.
Then $F$ belongs to $\D^{k,p}$ and $\norme{F}_{k,p} \leq C_{k,p} \limsup_n \norme{F_n}_{k,p}$ for some $C_{k,p} > 0$.
\end{lemma}

\begin{proof}[Proof of Lemma \ref{lem:Nualart}]
See Appendix \ref{appendixA}.
\end{proof}

We now have the ingredients needed to prove Proposition~\ref{prop:dinfty}. The proof of this result is rather standard, but needs a bit of care as the most convenient way of obtaining Malliavin smooth random variables is truncating the Karhunen--Lo\`eve expansion of $\Gamma$. Doing so we face the issue that there is no Onsager inequality available for this approximation of the field that we are aware of. We will bypass this difficulty by considering a further convolution of this truncated version of $\Gamma$ against a smooth mollifier $\varphi$ and then use the Onsager inequality \eqref{eq:lem_Onsager_conv} for convolution approximations.

\begin{proof}[Proof of Proposition~\ref{prop:dinfty}]
Here, we sketch the proof and give full details in the Appendix \ref{appendixB}. We start by showing that $M$ belongs to $\D^\infty$. 
Let $n \geq 1, \delta > 0, j \geq 0$ and $p \geq 1$.
In the following, we will denote
\[
\Gamma_\delta = \Gamma * \varphi_\delta, \quad \Gamma_{n,\delta} = \sum_{k=1}^n A_k e_k \ast \varphi_\delta,
\quad
M_\delta = \int_\C f(x) e^{i \beta \Gamma_\delta(x) + \frac{\beta^2}{2} \E[\Gamma_\delta(x)^2]} dx
\]
and
\[
M_{n,\delta} = \int_\C f(x) e^{i \beta \Gamma_{n,\delta}(x) + \frac{\beta^2}{2} \E[\Gamma_{n,\delta}(x)^2] } dx.
\]
$M_{n,\delta}$ is a smooth random variable (in the sense of Definition \ref{def:der}) and $D^j M_{n,\delta}$ is equal to
\[
(i \beta)^j \int_\C dx f(x) e^{i \beta \Gamma_{n,\delta}(x) + \frac{\beta^2}{2} \E[\Gamma_{n,\delta}(x)^2]} \sum_{k_1, \dots, k_j =1}^n (e_{k_1} \ast \varphi_\delta)(x) \dots (e_{k_j} \ast \varphi_\delta)(x) e_{k_1} \otimes \dots \otimes e_{k_j}.
\]
Combining Onsager inequalities, \eqref{eq:lem_bound_integral_1} and Lemma \ref{lem:Nualart}, one can show by taking the limit $n \to \infty$ that for all $k \geq 1$, $M_\delta \in \D^{k,2p}$ and that 
\[
\sup_{\delta > 0} \norme{M_\delta}_{k,2p} < \infty.
\]
Details of this are in the appendix. Now, because $(M_\delta, \delta >0)$ converges in $L^{2p}$ towards $M$, Lemma \ref{lem:Nualart} then implies that for all $k \geq 1$, $M \in \D^{k,2p}$. This concludes the proof that $M \in \D^\infty$. \\

The proof of the formula for $DM$ now follows via a series of approximation arguments. From the first part by taking $n \to \infty$, one can rather quickly deduce that
\[
DM_\delta = i \beta \int_\C dx f(x) e^{i\beta \Gamma_\delta(x) + \frac{\beta^2}{2} \E[\Gamma_\delta(x)^2]} \sum_{k=1}^\infty (e_k \ast \varphi_\delta)(x) e_k.
\]
Next, one argues that $(DM_\delta, \delta>0)$ converges in $L^2(\Omega;H)$ towards
\[
i \beta  \int_\C dx f(x) \mu(x) C(x, \cdot)
\]
and concludes that it necessarly corresponds to $DM$ by Lemma \ref{lem:Nualart0}. Here one again uses Onsager inequalities and dominated convergence. The full details are found in the appendix.
\end{proof}

\subsection{Bounds on $|\delta(DM)|$ and $\|D^2M\|_{\HC \otimes \HC}$ -- proof of Proposition \ref{prop:bound_delta_D2_intro}}\label{sec:bounds_for_derivatives}

The goal of this section is to control the tails of $|\delta(DM)|$ and $\|D^2M\|_{\HC \otimes \HC}$. We first note that these two random variables can be written explicitly in terms of imaginary chaos.

\begin{lemma}\label{lem:explicit_4_variables}
Let $f \in L^\infty(\C)$. Then
\begin{gather}
\label{eq:lem_explicit_a}
\delta(DM) = \beta \int_\C f(x) \frac{d}{d\beta} \mu(x) dx, \\
\label{eq:lem_explicit_c}
\|D^2M\|_{\HC \otimes \HC}^2 = \beta^4 \Re \int_{\C \times \C} f(x) f(y) \mu(x) \overline{\mu(y)} C(x,y)^2 dx dy,
\end{gather}
where the expression $\frac{d}{d\beta} \mu(x)$ is given sense by $\lim_{\delta \to 0}\left(i\Gamma_\delta(x)+\beta \E\Gamma_\delta^2(x)\right) :\exp(i\beta \Gamma_\delta(x)):$ with the limit, say, in $H^{-d}(U)$ and in probability.
\end{lemma}

The proof of \eqref{eq:lem_explicit_c} is very similar to the proof of the formula of $DM$ and we omit the details. The origin of  \eqref{eq:lem_explicit_a} can be explained by the following formal computation, that can be turned into a rigorous proof in a very similar manner as what we did in the proof of Proposition \ref{prop:dinfty} when we obtained the explicit expression of $DM$ -- one needs to use smooth approximations both for the field $\Gamma$, and smooth Malliavin variables.

\begin{proof}['Formal' proof of Lemma \ref{lem:explicit_4_variables}]

By Proposition \ref{prop:dinfty}, and then by integration by parts for $\delta$ (Proposition 1.3.3 of \cite{Nualart}), we have
\begin{align*}
\delta(DM) & = i \beta \int_\C f(x) \delta(\mu(x) C(x,\cdot)) dx \\
& = i \beta \int_\C f(x) \left( \mu(x) \delta(C(x,\cdot)) - \scalar{ D\mu(x), C(x,\cdot) }_{\HC} \right) dx.
\end{align*}
Noticing that $\delta(C(x,\cdot)) = \Gamma(x)$ (see (1.44) of \cite{Nualart}) and that by Proposition \ref{prop:dinfty}
$\scalar{ D\mu(x), C(x,\cdot) }_{\HC} = i \beta \mu(x) C(x,x)$, we obtain
\begin{align*}
\delta(DM) = \beta \int_\C f(x) \mu(x) (i \Gamma(x) + \beta C(x,x) ) dx = \beta \int_\C f(x) \frac{d}{d\beta} \mu(x) dx.
\end{align*}
This shows \eqref{eq:lem_explicit_a}.
\end{proof}

\begin{proof}[Proof of Proposition \ref{prop:bound_delta_D2_intro}]
We will only write the details for the variable $\delta(DM)$ since bounding the moments of $\|D^2M\|_{\HC \otimes \HC}$ is very similar to bounding the moments of imaginary chaos itself (with the use of \eqref{eq:lem_int2} instead of \eqref{eq:lem_bound_integral_1}). 

Let $N \ge 1$ and let $K \Subset U$ be the support of $f$. By Lemma~\ref{lem:explicit_4_variables} we have
\[\E[|\delta(DM)|^{2N}] \le \|f\|_\infty^{2N} \beta^{2N} \int_{K^{2N}} \Big| \E\Big[\prod_{j=1}^N \frac{d}{d\beta} \mu(x_j) \frac{d}{d\beta} \overline{\mu(y_j)}\Big] \Big| \, dx_1 \dots dx_N dy_1 \dots dy_N.\]
By a limiting argument, one can justify the formal identity:
\[\E\Big[\prod_{j=1}^N \frac{d}{d\beta} \mu(x_j) \frac{d}{d\beta} \overline{\mu(y_j)}\Big] = \Big[\prod_{\ell=1}^N \frac{d}{d\beta_\ell} \frac{d}{d\gamma_\ell} E((\beta_j)_{j=1}^N,(\gamma_j)_{j=1}^N)\Big]_{\beta_1 = \dots = \beta_N = \gamma_1 = \dots \gamma_N = \beta}.\]
where
\[E((\beta_j)_{j=1}^N,(\gamma_j)_{j=1}^N) \coloneqq e^{-\sum_{j < k} \beta_j \beta_k C(x_j,x_k) - \sum_{j < k} \gamma_j \gamma_k C(y_j,y_k) + \sum_{j,k} \beta_j \gamma_k C(x_j,y_k)}.\]
Let $(z_1,\dots,z_{2N}) \coloneqq (x_1,\dots,x_N,y_1,\dots,y_N)$.
By induction one sees that after differentiating w.r.t. the first $k$ of the variables $\beta_1,\dots,\beta_N,\gamma_1,\dots,\gamma_N$ and expanding one is left with a finite number of terms of the form
\[\pm \prod_{j=1}^{N} \beta_j^{n_j} \gamma_j^{m_j} \prod_{j=1}^\ell C(z_{a_j},z_{b_j}) E((\beta_j)_{j=1}^N,(\gamma_j)_{j=1}^N),\]
where $0 \le n_j,m_j,\ell \le k$, $1 \le a_1 < a_2 < \dots < a_\ell \le k$ and $1 \le b_1,\dots,b_\ell \le 2N$ with $a_j \neq b_j$ for all $j$.
Hence we have
\[\E[|\delta(DM)|^{2N}] \le C_N \sum_{\ell = 1}^{2N} \sum_{1 \le a_1 < \dots < a_\ell \le 2N}\sum_{b_1,\dots,b_\ell=1}^{2N} \int_{K^{2N}} \prod_{j=1}^\ell \1_{a_j \neq b_j} |C(z_{a_j}, z_{b_j})| e^{\mathcal{E}(\Gamma;\mathbf{x};\mathbf{y})} \, dx_j dy_j.\]
Note that $|C(z_{a_j}, z_{b_j})| \le C \log \frac{4R}{|z_{a_j} - z_{b_j}|}$ for some $C > 0$ and $R$ large enough so that $K \subset B(0,R)$. Thus applying Lemma~\ref{lem:Onsager_weak} to each summand, we can bound the whole sum by
\[C_N \int_{K^{2N}} \prod_{j=1}^{2N} \log \frac{4R}{\min_{k \neq j} |z_j - z_k|} (\min_{k \neq j} |z_j - z_k|)^{-\beta^2/2} \, dz_1 \dots dz_{2N}.\]
By scaling this is less than
\[C_N \int_{B(0,1/4)^{2N}} \prod_{j=1}^{2N} \log \frac{1}{\min_{k \neq j} |z_j - z_k|} (\min_{k \neq j} |z_j - z_k|)^{-\beta^2/2} \, dz_1 \dots dz_{2N},\]
which by Lemma~\ref{lem:integral_min} is less than $C_N (d - \beta^2)^{3N}$.
\end{proof}

\subsection{Small ball probabilities for the Malliavin determinant of $M$ -- proof of Proposition~\ref{prop:tail_det}}\label{sec:mainproof}

\noindent This section contains the main probabilistic input to Theorem~\ref{thm:main} -- the proof of Proposition~\ref{prop:tail_det}. 
Roughly, the content of this proposition is to establish super-polynomial decay of $\P(\det \gamma_M < \eps)$ as $\eps \to 0$, where 
$\det \gamma_M \coloneqq (\|D M\|_{\HC}^4 - |\langle DM, D\overline{M}\rangle_{\HC}|^2)/4$ is the Malliavin determinant of $M = \mu(f)$.

We will start by presenting a toy model explaining the strategy; then we explain the proof setup and prove the proposition modulo some technical chaos lemmas. The section finishes by proving the technical estimates.

\subsubsection{A toy model: small ball probabilities for $\|:\exp(i\beta \mathrm{GFF}):\|_{H^{-1}(\R^2)}$}\label{sec:toy1}

To explain the strategy of our proof, we consider a toy problem asking about the small ball probabilities for norms of imaginary chaos.
For concreteness, let us do it here with the 2D Gaussian free field; see Proposition \ref{prop:chaos_sobolev_tails} at the end of this section for a more general statement.
\begin{itemize}
    \item[] Consider the 2D zero boundary GFF on $K = [0,1]^2$ and the imaginary chaos $\mu_\beta$. We know that as a generalized function $\mu_\beta \in H^{-1}(K)$ for all $\beta \in (0, \sqrt{2})$. Can we prove super-polynomial bounds for $\P\left( \|\mu _\beta\|_{H^{-1}(K)} < \eps\right)$? Moreover, can we obtain bounds that are tight as $\beta \to \sqrt{2}$?
\end{itemize}

\noindent Writing out the norm squared, we have that
\[ \|\mu \|^2_{H^{-1}(K)} = \int_{K^2} \mu(x)G(x,y) \overline{\mu(y)}\,dx\,dy  > 0,\]
where $G$ is the Dirichlet Green's function on $K$. Now, the expectation $\E \|\mu \|^2_{H^{-1}(K)}$ is easy to calculate and it is bounded. As all moments exist, one could imagine proving bounds near zero by using concentration results on $\mu$. However, these concentration results do not see the special role of zero and would not suffice for good enough bounds for asymptotics near $0$. 

The idea is then to use only the decorrelated high-frequency part of $\Gamma$ to stay away from zero. To make this more precise, denote by $\Gamma_\delta$ the part of the GFF containing only frequencies less than $\delta^{-1}$ and let $\hat \Gamma_\delta = \Gamma - \Gamma_\delta$ denote the tail of the GFF. Consider now the projection bound $\|f\|_{H^{-1}(K)}\|\mu \|_{H^{-1}(K)} \geq \langle \mu, f \rangle_{H^{-1}(K)}$ for any $f \in H^{-1}(K)$. Setting $f(x) = f_\delta(x) = \Delta(:e^{i\beta \Gamma_\delta(x)}:)$, we get that
\[ \|\mu \|_{H^{-1}(K)} \geq \frac{|\int_{K} :e^{i\beta \hat \Gamma_\delta(x)}: e^{\beta^2 \E[\Gamma_{\delta}(x)^2]} \, dx| }{\|f_\delta\|_{H^{-1}(K)}}.\]
A small calculation shows that $\|f_\delta\|_{H^{-1}(K)} = \|:e^{i\beta \Gamma_\delta(y)}:\|_{H^1(K)}$. It is further believable that we should have $\|:e^{i\beta \Gamma_\delta(y)}:\|_{H^1(K)} \asymp \delta^{-\beta^2/2} \|\Gamma_\delta\|_{H^1(K)}$, and that this expression admits Gaussian concentration. As in the concrete case $\E \|\Gamma_\delta\|_{H^1(K)}  \asymp \delta^{-1}$, we can conclude that the denominator is of order $\delta^{-1-\beta^2/2}$ with super-polynomial concentration on fluctuations.

In the numerator, the term of the form $\int_K :e^{i\beta \hat \Gamma_\delta(x)}: e^{\beta^2 \E[\Gamma_{\delta}(x)^2]} dx $ remains. Such a tail chaos is very highly concentrated around its mean which is of order $\delta^{-\beta^2}$, with fluctuations of unit order having a super-polynomial cost in $\delta$. Thus the whole ratio will concentrate around 
\[C\frac{\delta^{-\beta^2}}{\delta^{-1-\beta^2/2}} \sim C\delta^{1-\beta^2/2},\]
with super-polynomial cost for fluctuations on the same scale. Thus setting $\eps = \delta^{1-\beta^2/2}$ we obtain super-polynomial decay for $\P\left( \|\mu \|_{H^{-1}(K)} < \eps\right)$. 

Whereas this is good enough for any fixed $\beta$, observe that as $\beta \to \sqrt{2}$ the exponent $1 - \beta^2/2$ goes to $0$. Moreover, we have $\E \|\mu \|_{H^{-1}(K)}^2 = O((2-\beta^2)^{-2})$, but $\E |\int :e^{-i\beta \hat \Gamma_\delta(x)}: |^2 = O((2-\beta^2)^{-1})$. As further $\|f_\delta\|_{H^{-1}(K)} \asymp \delta^{-\beta^2/2}\|\Gamma_\delta\|_{H^1(K)}$ and $\|\Gamma_\delta\|_{H^1(K)}$ does not depend on $\beta$, we see that we are in fact losing in terms of $\beta^2-2$ as well. 

Illustratively, we are losing in high frequencies because we are replacing 
\[\int \mu(x)G(x,y) \overline{\mu(y)}\quad\quad\quad\text{by}\quad\quad\quad\int :e^{i\beta \hat \Gamma_\delta(x)}: :e^{-i\beta \hat \Gamma_\delta(y)}:.\] 
After taking expectation, in terms of near-diagonal contributions, as $G(x,y) \sim -\log |x-y|$ near the diagonal, this basically translates to replacing $-\int |x|^{-\beta^2/2} \log |x|$ with $\int |x|^{-\beta^2/2}$ and results in the loss of a factor of $2-\beta^2$ as $\beta^2 \to 2$. Thus we have to tweak our test function $f_\delta$ further to at the same time guarantee sufficient concentration and not to lose too much on tails.

We will see later on that this strategy gives us more generally the following result.

\begin{proposition}\label{prop:chaos_sobolev_tails}
  Let $f \in C_c^\infty(U)$. Then for each $\nu \in (0, \sqrt{d})$, there exist constants $c_1,c_2,c_3 > 0$ such that
  \[\P[\|f \mu\|_{H^{-d/2}(\reals^d)} \le (d-\beta^2)^{-2} \lambda] \le c_1 e^{-c_2 \lambda^{-c_3}}\]
  for all $\lambda > 0$ and all $\beta \in (\nu, \sqrt{d})$.
\end{proposition}

The same strategy for the determinant requires some extra input, yet the key ideas are present already in this toy model: the projection bound corresponds to the analogue of Malliavin determinants given by Lemma \ref{lem:deltalowerbound}, the concentration of the numerator to Lemma \ref{lemma:mainterm} and that of the denominator to Lemma \ref{lemma:normterm}. The only new technical ingredient will enter as Lemma \ref{lemma:oddterm}.
  \\

\subsubsection{Proof setup and proof of Proposition \ref{prop:tail_det} modulo technical lemmas}

Let $f$ be a bounded continuous function whose support is a compact subset of $U$ and set $M = \mu(f)$.
Our goal in this section is to obtain lower bounds on $\P[\det \gamma_M \ge \lambda]$, where $\det \gamma_M$ is the Malliavin determinant \eqref{eq:malliavin_det}.

As in the toy problem, it is not so clear how to obtain sharp bounds directly and the idea is to use the projection bound from Lemma~\ref{lem:deltalowerbound}, which says that
\begin{equation}\label{eq:P_after_projection}\P[\det \gamma_M \ge \lambda] \ge \P\Big[\frac{(|\langle DM, h\rangle_{H_\C}| - |\langle D\overline{M}, h\rangle_{H_\C}|)^4}{\|h\|_{H_\C}^4} \ge 4 \lambda\Big]\end{equation}
for any $h \in H_\C$.
A key step is the specific choice of $h(x)$, which needs to at the same time give a precise enough bound and allow for chaos computations. Moreover, we have to ensure that it also belongs to the Cameron--Martin space. Here, one of the technical difficulties is that in general we do not have a good understanding of the Cameron--Martin space of $\Gamma$. To deal with that, we will use the decomposition theorem, Theorem \ref{thm:decomposition} to be able to work with almost $\star$-scale invariant fields.

More precisely, let us fix an open set $V$ with $\overline{V}$ a compact subset of $U$ such that $\supp f \subset V$. Then by Theorem~\ref{thm:decomposition} one can write $\Gamma|_V = Y + Z \eqqcolon X$ where $Y$ is an almost $\star$-scale invariant field with smooth and compactly supported seed covariance $k$ and parameter $\alpha$, and $Z$ is an independent Hölder-continuous field. Recall further the approximations $Y_\eps$ of $Y$ of such a field from Section \ref{sec:starscale} and the notation for its tail field $\hat{Y}_\varepsilon \coloneqq Y - Y_\varepsilon$. 

Now, notice that
\[\det \gamma_M = \frac{\beta^4}{4} \Big(\Big|\int_U f(x)f(y) \mu(x) \overline{\mu(y)} C(x,y) \, dx \, dy\Big|^2 - \Big|\int_U f(x)f(y) \mu(x) \mu(y) C(x,y) \, dx \, dy\Big|^2\Big),\]
where the right hand side only depends on $\mu$, and thus on $\Gamma$, restricted to $V$. Thus, to obtain bounds on $\det \gamma_M$, we can instead of working with the (complexified) Cameron--Martin space $H_\C = H_{\Gamma, \C}$, just as well work with the Cameron--Martin space of $Y + Z$, which is defined on the whole plane. Apologising for the abuse of notation, we still denote it by $H_\C$. This small trick allows us to use the independence structure of the field $Y$, and also puts Fourier techniques in our hand.

\paragraph*{Definition of $h$.}

Whereas the decomposition theorem and the change of Cameron--Martin space make the computations potentially doable, they become practically doable only with a very careful choice of the test function $h$. Namely, we set
\[h(x) = h_{\delta}(x) = e^{i\beta Y_\delta(x) - \frac{\beta^2}{2} \E[Y_\delta(x)^2]} \int_U f(y) :e^{i\beta Z(y)}: :e^{i\beta\hat{Y}_\delta(y)}: R_\delta(x,y) \, dy,\]
where $R_\delta(x,y) = g_\delta(x)g_\delta(y)\E[\hat{Y}_\delta(x)\hat{Y}_\delta(y)]$ is defined using a smooth indicator $g_\delta$ of $\delta$-separated squares and the parameter $\delta$ will be chosen in a suitable way according to $\lambda$. 

More precisely, let $\Qc_\delta$ be the collection of cubes of the form \[[4k_1\delta,(4k_1+1)\delta] \times \dots \times [4k_d\delta, (4k_d+1)\delta],\]
where $k_1,\dots,k_d \in \integers$.
Note in particular that the cubes are $\delta$-separated and hence the restrictions of $\hat{Y}_\delta$ to two distinct cubes in $\Qc_\delta$ are independent.
We then set 
\begin{equation}\label{eq:g_delta}
g_\delta = \varphi_\delta * \1_{\bigcup \Qc_\delta \cap V},    
\end{equation}
where $\varphi$ is a smooth mollifier supported in the unit ball and $\varphi_\delta(x) = \delta^{-d} \varphi(x/\delta)$.

We note that $h$ is indeed almost surely an element of ${H_\C}$, since the Malliavin derivative of $(i\beta)^{-1}\int f(y) :e^{i\beta Z(y)}: g_\delta(y) :e^{i\beta \hat{Y}_\delta(y)}: \, dy$ with respect to the field $\hat{Y}_\delta$ equals $$x \mapsto \int_U f(y) :e^{i\beta Z(y)}: g_\delta(y) :e^{i\beta \hat{Y}_\delta(y)}: \E[\hat{Y}_\delta(x) \hat{Y}_\delta(y)] \, dy$$
and lies in $H_{\hat{Y}_\delta,\C}$ (the complexification of the Cameron--Martin space of $\hat{Y}_\delta$). In particular, since $Y = Y_\delta + \hat{Y}_\delta$ is an independent sum, it lies in $H_{Y,\C}$ as well and, by Lemma~\ref{lemma:C_Y_fourier}, this as a set of functions coincides with $H_\C^{d/2}(\reals^d)$.
Moreover, the map $x \mapsto g_\delta(x) e^{i\beta Y_\delta(x) - \frac{\beta^2}{2} \E[Y_\delta(x)^2]}$ is almost surely smooth so multiplying by it shows that
\[x \mapsto g_\delta(x) e^{i\beta Y_\delta(x) - \frac{\beta^2}{2} \E[Y_\delta(x)^2]} \int_U f(y) :e^{i\beta Z(y)}: g_\delta(y) :e^{i\beta \hat{Y}_\delta(y)}: \E[\hat{Y}_\delta(x) \hat{Y}_\delta(y)] \, dy \in H^{d/2}_\C(\reals^d).\]
Finally, as $Y + Z$ is an independent sum, Lemma \ref{lem:CMdecomp} implies that $H_\C^{d/2}(\reals^d) \subset H_{\C}$ as desired. 

\paragraph*{Proof of Proposition \ref{prop:tail_det}}

In order to derive bounds on $\P[\det \gamma_M < \lambda]$ and $\P(\frac{\det \gamma_M}{\|DM\|^2_{\HC}} < \lambda)$ for $\lambda > 0$ small, we will look at the three terms $|\langle DM, h_\delta\rangle_{\HC}|$, $|\langle D\overline{M}, h_\delta\rangle_{\HC}|$ and $\|h_\delta\|_{\HC}$ appearing in \eqref{eq:P_after_projection} separately and collect the results in the following lemmas. 

  \begin{lemma}\label{lemma:mainterm}
	For every $\nu >0$, there exists a constant $c_2 > 0$ such that for all $c > 0$ small enough $$\P[|\langle DM, h_\delta \rangle_{H_\C}| \le c (d-\beta^2)^{-2}\delta^{d}] \le \exp \left( -c_2\delta^{-{d\wedge 2}} \right)$$
	for all small enough $\delta > 0$ and all $\beta \in (\nu, \sqrt{d})$.
\end{lemma}

\begin{lemma}\label{lemma:normterm}
For all $\eta > 0$ small enough, we can choose $C > 0$ such that $$\|h_\delta\|^2_{H_\C} \leq C\delta^{\beta^2-2d-2\eta}W^2 |\langle DM, h_\delta \rangle_{H_\C}|,$$ where $W$ is a $Y_\delta$-measurable positive random variable. Moreover, we can pick $c_1, c_2 > 0$ such that for all  $\delta \in (0,1)$ and $t \ge c_1\delta^{-2-\eta}$ we have
\[\P(W > t) \leq \exp(-c_2\delta^{\eta}t^{\frac{2}{d}}).\]
\end{lemma}

\begin{lemma}\label{lemma:oddterm}
For every $\nu >0$, there exists a constant $c_1 > 0$ such that the following holds. For every $c > 0$, we can choose $c_2 > 0$ such that $$\P[|\langle \overline{DM}, h_\delta \rangle_{H_\C}| \ge c (d-\beta^2)^{-2}\delta^{d}] \le \exp(-c_2\delta^{-c_1})$$
for all small enough $\delta > 0$ and all $\beta \in (\nu, \sqrt{d})$.
\end{lemma}

We now explain how we deduce Proposition \ref{prop:tail_det} from these lemmas, and then in the next subsections turn to their proofs.

\begin{proof}[Proof of Proposition \ref{prop:tail_det}]

	By Lemma \ref{lem:deltalowerbound}, we have that 
		\[\P\left(\frac{\det \gamma_M}{\|DM\|_{\HC}^2} \geq \eps / 4\right) \geq \P \left( \frac{(|\langle DM, h_\delta \rangle_{\HC}| - |\langle D\overline{M}, h_\delta \rangle_{\HC}|)^2}{\|h_\delta\|_{\HC}^2}\ge \varepsilon \right)\]
	and
	\[\P \left( \det \gamma_M \geq \eps / 4 \right) \geq \P \left( \frac{(|\langle DM, h_\delta \rangle_{\HC}| - |\langle D\overline{M}, h_\delta \rangle_{\HC}|)^2}{\|h_\delta\|_{\HC}^2} \ge \sqrt{\varepsilon} \right),\]
	so it suffices to bound $\P \left( \frac{(|\langle DM, h_\delta \rangle_{\HC}| - |\langle D\overline{M}, h_\delta \rangle_{\HC}|)^2}{\|h_\delta\|_{\HC}^2}\le \varepsilon \right)$ from above. Here $h_\delta$ is as above and we will choose  $\delta$ depending on $\eps$.
   
   Using Lemma \ref{lemma:normterm}, we first bound for some $\eta > 0$
\[\frac{(|\langle DM, h_\delta \rangle_{\HC}| - |\langle D\overline{M}, h_\delta \rangle_{\HC}|)^2}{\|h_\delta\|_{\HC}^2} \geq C^{-1} \delta^{-\beta^2+2d+2\eta}W^{-2}(|\langle DM, h_\delta \rangle_{\HC}|-2|\langle D\overline{M}, h_\delta \rangle_{\HC}|).\]
   Hence, taking $c$ to be the constant from Lemma \ref{lemma:mainterm} we can bound
   \[\P\Big( \frac{(|\langle DM, h_\delta \rangle_{\HC}| - |\langle D\overline{M}, h_\delta \rangle_{\HC}|)^2}{\|h_\delta\|_{\HC}^2}\le  (d-\beta^2)^{-2}\delta^{3d+5}\Big)\]
   by
   \[\P \left(|\langle DM, h_\delta \rangle_{\HC}|-2|\langle D\overline{M}, h_\delta \rangle_{\HC}| \leq  \frac{c}{2}(d-\beta^2)^{-2}\delta^{d} \right) +
   \P \left( C\delta^{\beta^2-2d-2\eta}W^2 >  \frac{c}{2}\delta^{-2d-5} \right).\] 
   The second term can be bounded using Lemma \ref{lemma:normterm} loosely
   by $\exp(-c_1\delta^{-c_1})$ for some $c_1 > 0$.
   
   For the first term, Lemma \ref{lemma:mainterm} gives that 
   $$\P(|\langle DM, h_\delta \rangle_{\HC}| \leq  c(d-\beta^2)^{-2}\delta^{d}) \leq \exp(-c_2\delta^{-d\wedge 2})$$
   and Lemma \ref{lemma:oddterm} gives constants $c_3 > 0$
   $$\P(2|\langle \overline{DM}, h_\delta \rangle_{\HC}| \geq  \frac{c}{2}(d-\beta^2)^{-2}\delta^{d}) \leq \exp(-\delta^{-c_3}),$$
 and we thus obtain the proposition.
	
The case of the standard log-correlated field on circle needs extra attention, and is treated in Section \ref{sec:circlecase}.
\end{proof}

One can see that a simplified version of the above proof can also be used to prove Proposition~\ref{prop:chaos_sobolev_tails}.

\begin{proof}[Proof of Proposition~\ref{prop:chaos_sobolev_tails}]
Recall that on the support of $f$, we can write $\Gamma_{|V} = Y + Z = X$, where $Y$ is almost $\star-$scale invariant and $Z$ is Holder regular, both defined on the whole space.
Note that by Lemma~\ref{lemma:C_Y_fourier} and Theorem~\ref{thm:decomposition} the operators $C_Y$ and $C_Z$ are bounded from $H^{-d/2}(\reals^d)$ to $H^{d/2}(\reals^d)$ and hence so is $C_X$.
Thus for any $\varphi \in H^{-d/2}(\reals^d)$ we have
\[\langle C_X \varphi, \varphi \rangle_{L^2(\reals^d)} \le \|C_X \varphi\|_{H^{d/2}(\reals^d)} \|\varphi\|_{H^{-d/2}(\reals^d)} \le \|C_X\|_{H^{-d/2}(\reals^d) \to H^{d/2}(\reals^d)} \|\varphi\|_{H^{-d/2}(\reals^d)}^2\]
so that in particular
\[\|f\mu\|_{H^{-d/2}(\reals^d)}^2 \gtrsim \langle C_X (f \mu), f\mu \rangle_{L^2(\reals^d)} = \beta^{-2} \|DM\|_{\HC}^2 \ge \beta^{-2} \frac{|\langle DM, h_\delta\rangle_{\HC}|^2}{\|h_\delta\|_{\HC}^2}.\]
Using this inequality one can proceed as in the proof of Proposition~\ref{prop:tail_det} except one does not need to take care of the term $\langle D\overline{M}, h_\delta\rangle$.
\end{proof}

The rest of this subsection is dedicated to the proofs of Lemmas \ref{lemma:mainterm}, \ref{lemma:normterm} and \ref{lemma:oddterm}, and sketching the extension to the case of the circle. 

\subsubsection{Proof of Lemma~\ref{lemma:mainterm}}

\begin{proof}[Proof of Lemma~\ref{lemma:mainterm}]
  Let us fix some $\nu > 0$  small.
  Note that $\langle DM, h_\delta \rangle_{\HC}$ is equal to
  \begin{align*}
    & i\beta \int_U f(x) :e^{i\beta X(x)}: \overline{h_\delta(x)} \, dx = i\beta \int_{U\times U} f(x) f(y) :e^{i\beta (\hat{Y}_\delta(x) + Z(x))}: :e^{-i\beta (\hat{Y}_\delta(y) + Z(y))}: R_\delta(x,y) \\
    & = i\beta\sum_{Q \in \Qc_\delta} \int_{Q \times Q} f(x)f(y) :e^{i\beta(\hat{Y}_\delta(x) + Z(x))}: :e^{-i\beta(\hat{Y}_\delta(y) + Z(y))}: R_\delta(x,y) \, dx \, dy
  \end{align*}
  since $R_\delta(x,y) = 0$ if $x$ and $y$ are not in the same square in $\Qc_\delta$.  Moreover the summands are mutually independent,
  when we condition on the field $Z$, and by scaling each term agrees in law with
  \[\delta^{2d} J_Q \coloneqq \delta^{2d} \int_{\delta^{-1}Q \times \delta^{-1}Q} f(\delta x) :e^{i\beta Z(\delta x)}: :e^{-i\beta Z(\delta y)}:f(\delta y) :e^{i\beta \hat{Y}_\delta(\delta x)}: :e^{-i\beta \hat{Y}_\delta(\delta y)}: R_\delta(\delta x,\delta y) \, dx \, dy.\]
We can write 
  \[\E[J_Q|Z] = \int_{\delta^{-1} Q \times \delta^{-1} Q} f(\delta x) f(\delta y) :e^{i\beta Z(\delta x)}: :e^{-i\beta Z(\delta y)}: e^{\beta^2 \E[\hat{Y}_\delta(\delta x)\hat{Y}_\delta(\delta y)]} R_\delta(\delta x,\delta y) \, dx \, dy.\]
  Whenever $Q$ is such that $f(x) \ge \|f\|_\infty/2$ for all $x \in Q$ (or similarly if $f(x) \le - \|f\|_\infty/2$), and the event $E_Q := \{\sup_{x,y \in Q} |Z(x)-Z(y)| \leq \pi/(4\beta)\}$ holds, a basic calculation that uses Lemma \ref{lemma:star_scale_cov_2} shows that
  \begin{itemize}
      \item $\E[J_Q| Z, E_Q] \ge C(d-\beta^2)^{-2}$,  for some constant $C > 0$ that is uniform over $\beta \in (\nu, d)$ and depends only on $\|f\|_\infty$
      \item $\E[J_Q^2| Z, E_Q] \le c(d-\beta^2)^{-4}$ for some constant $c > 0$ that is again uniform over $\beta \in (\nu, d)$ and depends solely on $\|f\|_\infty$. 
      \end{itemize}
In particular, by the Paley-Zygmund inequality for any such square $Q$  it holds that $\P[J_Q \ge \lambda(d-\beta^2)^{-2}|Z, E_Q] \ge p$, where $\lambda = C/2$ and $p > 0$ is some constant. In the following, we denote by $\tilde{\Qc}_\delta$ the collection of those squares in which $f$ is larger than $\|f\|_\infty/2$ (again, we may consider $-f$ instead of $f$ if needed).
  
  Now, recall that $Z$ is a Hölder continuous Gaussian field, and thus by local chaining inequalities (e.g. Proposition 5.35 in \cite{VanHandel}), we have that for some universal constant $C > 0$
  $$\P \left(\sup_{|x-y|\leq 2\delta}|Z(x)-Z(y)| > \pi/(4\beta) \right) \leq C\exp(-C\delta^{-2}).$$  Thus denoting $E = \{\sup_{|x-y|\leq 2\delta}|Z(x)-Z(y)| \leq \pi/(4\beta)\}$ , we can bound 
  $$\P[|\langle DM, h_\delta \rangle_{\HC}| \le c (d-\beta^2)^{-2}\delta^{d}] \leq P(E^c) + \P\Big[|\langle DM, h_\delta \rangle_{\HC}| \le c (d-\beta^2)^{-2}\delta^{d}|E\Big].$$
 As $\P(E^c) \leq C\exp(-C\delta^{-2})$ and $E \subseteq \bigcap_Q E_Q$, it remains to only take care of the second term working under the assumption that the event $E_Q$ holds for all $Q$. For any $t>0$ to be chosen later, we have
  \begin{align*}
    \P\Big[|\langle DM, h_\delta\rangle_H| \le (d-\beta^2)^{-2}t|E\Big] & \le \P\Big[J_Q \ge (d-\beta^2)^{-2}\lambda \text{ for at most } t/(\beta \lambda \delta^{2d}) \text{ distinct } Q \in \tilde{\Qc}_\delta|E\Big] \\
    & \le \P[\Bin(|\tilde{\Qc}_\delta|, p) \le t/(\beta \lambda \delta^{2d})] \\
    & \le e^{-2|\tilde{\Qc}_\delta| \Big(p - \left\lceil\frac{t}{\beta \lambda \delta^{2d}} \right\rceil |\tilde{\Qc}_\delta|^{-1} \Big)^2}
  \end{align*}
  where $\Bin(n,p)$ denotes the Binomial distribution. In the second line we used the conditional independence of $J_Q$ given $Z$ and the conditional probability obtained above; on the last line we used the Hoeffding's inequality
  \[\P[\Bin(n,p) \le m] \le e^{-2n(p - \frac{m}{n})^2}.\]
  Noting that $c_1 \delta^{-d} \le |\tilde{\Qc}_\delta| \le c_2 \delta^{-d}$ for some $c_1,c_2 > 0$, we see that by choosing $t = p \beta \lambda \delta^d / (2 c_2)$ we get 
  \[\P\Big[|\langle DM, h_\delta\rangle_H| \le (d-\beta^2)^{-2}t|E \Big] \le e^{-2 c_1 \frac{p}{3} \delta^{-d}}\] 
  for small enough $\delta > 0$ and the lemma follows.
\end{proof}

\subsubsection{Proof of Lemma~\ref{lemma:normterm}}

\begin{proof}[Proof of Lemma~\ref{lemma:normterm}]
	We start with some immediate bounds that allow the usage of inequalities on Sobolev spaces $H_\C^s(\R^d)$. First, by Lemma~\ref{lemma:C_Y_fourier} we have
			\[C^{-1}\| \cdot \|_{H^{d/2}_\C(\R^d)} \leq \| \cdot\|_{H_{Y,\C}} \leq C\| \cdot \|_{H^{d/2}_\C(\R^d)}\]
	for some $C > 0$.
	On the other hand, by Lemma \ref{lem:CMdecomp}, we have that
	$$\| \cdot \|_{\HC} \leq \| \cdot \|_{H_{Y,\C}} \leq \| \cdot \|_{H_{\hat Y_\delta,\C}}.$$
	Now let $\psi \in C_c^\infty(\reals^d)$ be a non-negative function which equals $1$ in the support of $g_\delta$ (recall that $g_\delta$ is defined in \eqref{eq:g_delta}). Set
	\[F(x) \coloneqq e^{i\beta Y_\delta(x)- \frac{\beta^2}{2} \E[Y_\delta(x)^2]} \psi(x)\]
	and 
	\[G(x) \coloneqq \int_{U \times U} f(y) :e^{i\beta Z(y)}: :e^{i\beta \hat{Y}_\delta(y)}: g_\delta(y) \E[\hat{Y}_\delta(x) \hat{Y}_\delta(y)] \, dx \, dy\]
	so that $g_\delta(x)F(x)G(x) = h_\delta(x)$.
    Using the above norm bounds in conjunction with the classical inequality $\|FG\|_{H^{d/2}(\R^d)} \lesssim \|F\|_{H^{d/2+\eps}_\C(\R^d)}\|G\|_{H^{d/2}_\C(\R^d)}$ for any $\eps > 0$ (see e.g. Theorem 5.1 in \cite{BH}), we can bound $\|h_\delta\|_{\HC}$ by some constant times
    	\[\|gFG\|_{H^{d/2}_\C(\R^d)} \lesssim \|g_\delta\|_{H^{d/2+\eps}_\C(\R^d)}\|F\|_{H^{d/2+\eps}_\C(\R^d)}\|G\|_{H^{d/2}_\C(\R^d)} \lesssim \|g_\delta\|_{H^{d/2+\eps}_\C(\R^d)}\|F\|_{H^{d/2+\eps}_\C(\R^d)}\|G\|_{H_{\hat Y_\delta,\C}}.\]
    We can bound $\|g_\delta\|_{H^{d/2+\eps}_\C(\R^d)} \lesssim \delta^{-d-\eps}$ by scaling and triangle inequality. Further, by definition we have that $\|G\|^2_{H_{\hat Y_\delta,\C}}= |\langle DM, h_\delta \rangle_{\HC}|$. Thus it remains to deal with $\|F\|_{H^{d/2+\eps}_\C(\R^d)}$. To do this, we will use Gaussian concentration inequalities. 
    
    Namely, by Theorem 4.5.7 in \cite{BogachevGauss}, if $X$ is isonormal on a Hilbert space $H'$, and any $T: H' \to \R$ is $L-$Lipschitz w.r.t $\| \cdot \|_{H'}$, then for all $t > 0$
	\[\P(T(X) - \E T(X) > t) \leq \exp(-\frac{t^2}{2L^2}).\]
We will make use of this concentration in the case $T = \|\cdot\|_{H^{d/2+\varepsilon}(\reals^d)}$ to bound $W := T(F)$.
We first apply Theorem A in \cite{AdamsFrazier}, which gives that for $f \in H^{d/2 + \eps}(\reals^d)$ we have $\|\exp(i f)\psi\|_{H^{d/2 + \eps}_\C} \lesssim \|f\|_{H^{d/2 + \eps}(\reals^d)} + \|f\|_{H^{d/2 + \eps}(\reals^d)}^{d/2 + \eps}$.\footnote{In \cite{AdamsFrazier} the authors consider compositions with real-valued functions; in our case one can apply it directly to the real and imaginary part. Note that by the theorem the first operator in the chain $f \mapsto e^{i f} - 1 \mapsto (e^{i f} - 1)\psi \mapsto e^{i f} \psi$ is bounded and the other two are bounded since $\psi$ is smooth.} This together with the fact that $\E[Y_{\delta}(x)^2]$ is constant in $x$ gives us that
$\|F\|_{H^{d/2+\eps}_\C(\R^d)} \leq c\delta^{\beta^2/2}(\|Y_\delta \tilde{\psi}\|_{H^{d/2+\eps}(\reals^d)} + \|Y_\delta \tilde{\psi}\|_{H^{d/2 + \eps}(\reals^d)}^{d/2 + \eps})$
for some $c > 0$. Here $\tilde{\psi} \in C_c^\infty(\reals^d)$ is some function which is $1$ in the support of $\psi$. Further, we have the following bounds:

\begin{claim}\label{claim1}
It holds that

1. $\| \cdot \|_{H^{d/2+\eps}(\reals^d)}$ is $O(\delta^{-2\eps})-$Lipschitz with respect to $\| \cdot \|_{H_{Y_\delta}}$.

2. $(\E \|\tilde{\psi} Y_\delta\|_{H^{d/2+\eps}(\reals^d)})^2 \leq \E \|\tilde{\psi} Y_\delta\|^2_{H^{d/2+\eps}(\reals^d)} \lesssim \delta^{-d-4\eps}$.
\end{claim}

\begin{proof}[Proof of Claim \ref{claim1}]
Recall from the proof of Lemma~\ref{lemma:C_Y_fourier} that the operator $C_{Y_\delta}$ is a Fourier multiplier operator with the symbol
\[
\hat{K}_\delta(\xi) := \int_\delta^1 v^{d-1} (1-v^\alpha) \hat{k}(v \xi) dv
\]
and $k$ is by assumption smooth. Moreover,
\[\|f\|_{H_{Y_\delta}}^2 = \int_{\reals^d} \hat{K}_\delta(\xi)^{-1} |\hat{f}(\xi)|^2 \, d\xi\]
and
\[\E\|\tilde{\psi} Y_\delta\|_{H^{d/2 + \varepsilon}(\reals^d)}^2 = \int_{\reals^d} (1 + |\xi|^2)^{d/2 + \varepsilon} \int_{\reals^d} |\hat{\tilde{\psi}}(\zeta)|^2 \hat{K}_\delta(\xi - \zeta) \, d\zeta \, d\xi.\]
The two claims thus directly follow from bounding $\hat{K}_\delta$ respectively by
\begin{align}
\label{eq:chat1}
    \hat{K}_\delta(\xi) & \lesssim \delta^{-2\eps} (1 + |\xi|^2 )^{-d/2 - \eps},\\
    \text{and} \quad \hat{K}_\delta(\xi) & \lesssim \delta^{-d-4\eps} (1 + |\xi|^2 )^{-d - 2\eps},
    \label{eq:chat2}
\end{align}
where the underlying constants do not depend on $\delta$. These inequalities are clear when $|\xi| \leq 1$, and follow by integrating the bounds  $\hat{k}(v \xi) \leq C |v \xi|^{-d-2\eps}$ and $\hat{k}(v \xi) \leq C |v \xi|^{-2d-4\eps}$ for $|\xi| > 1$.
\end{proof}

We can finally apply the Gaussian concentration to deduce that for all $\eps \in (0,d/2)$, there are some $c,C' > 0$, such that for all $t > c\delta^{-d-4\eps}$
	\[\P(\|\tilde{\psi} Y_\delta\|_{H^{d/2+\eps}(\R^d)} > t) \leq \exp\left(-C'\delta^{\eps}t^2\right),\]
	and thus for some $c',C'' > 0$ and for all $t > c'\delta^{-2-4\eps}$
\[\P(\|\tilde{\psi} Y_\delta\|_{H^{d/2+\eps}(\R^d)}+\|\tilde{\psi} Y_\delta\|_{H^{d/2+\eps}(\R^d)}^{d/2+\eps} > t) \leq \exp\left(-C'\delta^{\eps}t^{\frac{2}{d}}\right),\]
	implying the lemma. 
\end{proof}

\subsubsection{Proof of Lemma~\ref{lemma:oddterm}}

\begin{proof}
	We have
	\begin{align*}
	\scalar{DM, \overline{h_\delta}} = i \beta \int_{U \times U} f(x) f(y) e^{-2\beta^2 \E[X_\delta(x)^2]} : e^{i 2\beta X_\delta(x)} : :e^{i \beta \hat{Y}_\delta(x)}: :e^{i \beta \hat{Y}_\delta(y)}: R_\delta(x,y) dx dy,
	\end{align*}
	which we can write as a sum
	\[i \beta \sum_{Q \in \Qc_\delta} \int_{Q\times Q} f(x) f(y) e^{-2\beta^2 \E[X_\delta(x)^2]} : e^{i 2\beta X_\delta(x)} : :e^{i \beta \hat{Y}_\delta(x)}: :e^{i \beta \hat{Y}_\delta(y)}: R_\delta(x,y) dx dy =: i \beta \sum_{Q \in \Qc_\delta} L_Q.\]
	We can then first bound $$\E|\scalar{\overline{DM}, h_\delta}|^{2N} \le \beta^{2N} \E |\sum_{Q \in \Qc_\delta}L_Q|^{2N}.$$ If we expand the $2N$-th moment of such a sum, we obtain terms of the form
	\[\beta^{2N} \E\Big[ L_{Q_1} \dots L_{Q_N} \overline{L_{Q'_1} \dots L_{Q'_N}}\Big].\]
	Before taking expectation in each such term we separate the field $Y_\delta = Y_{\sqrt{\delta}} + \widetilde Y_\delta$, with $\widetilde Y_\delta := Y_\delta - Y_{\sqrt{\delta}}$ being independent of $Y_{\sqrt{\delta}}$. 
	We can then write each term as
	\begin{align*}
	& = \beta^{2N} \int_{U^{2N}} \prod_{j=1}^N f(x_j) f(y_j) f(x_j') f(y_j') R_\delta(x_j,y_j) R_\delta(x_j',y_j') e^{4\beta^2 \mathcal{E}(Y_{\sqrt{\delta}};\mathbf{x};\mathbf{x'})} e^{\beta^2 \mathcal{E}(\hat{Y}_{\delta};\mathbf{x},\mathbf{y};\mathbf{x}',\mathbf{y}')} \\
	& \quad \times e^{-2\beta^2 \sum_{j=1}^N (\E[X_\delta(x_j)^2] + \E[X_\delta(x_j')^2])} \E\left(\prod_{j=1}^{N} :e^{i2\beta (Z(x_j)+ \widetilde Y_\delta(x_j))}: :e^{-i2\beta(Z(x_j') + \widetilde Y_\delta(x_j'))}:\right),
	\end{align*}
	where the integration is over $x_j,y_j \in Q_j$ and $x_j',y_j' \in Q_j'$.
	We bound the expectation by
	\[\E \left|\prod_{j=1}^{N} :e^{i2\beta (Z(x_j)+ \widetilde Y_\delta(x_j))}: :e^{-i2\beta(Z(x_j') + \widetilde Y_\delta(x_j'))}:\right| \le C^N \delta^{-2N\beta^2},\]
	since $\E[\widetilde Y_\delta(x)^2] = \frac{1}{2} \log \frac{1}{\delta} + O(1)$.
	Now, there is some $c > 0$ such that $\mathcal{E}( Y_{\delta^{1/2}}; \mathbf{x} ; \mathbf{x}') \geq \mathcal{E}(Y_\delta^{1/2}, \mathbf{q}; \mathbf{q}') - c \sqrt{\delta} N^2$, where $\mathbf{q}$ and $\mathbf{q}'$ denote the vectors of midpoints for the ordered squares $Q_j$ and $Q_j'$. This can be seen by noting that since the seed covariance $k$ is Lipschitz, we have
	\[|\E[Y_{\sqrt{\delta}}(x)Y_{\sqrt{\delta}}(x')] - \E[Y_{\sqrt{\delta}}(q)Y_{\sqrt{\delta}}(q')]| \lesssim \int_0^{\frac{1}{2}\log \frac{1}{\delta}} e^u | |x-x'| - |q - q'| | (1 - e^{-\alpha u}) \, du \lesssim \sqrt{\delta}\]
	when $|x-q|, |x'-q'| \lesssim \delta$.
	Thus we obtain the upper bound
	$$\|f\|_\infty^{4N} \beta^{2N}\delta^{2\beta^2N}e^{c \sqrt{\delta} N^2} e^{4 \beta^2 \mathcal{E}( Y_{\delta^{1/2}}; \mathbf{q}_1 ; \mathbf{q}_2 )}\E[J_{Q_1} \dots J_{Q_N} \overline{J_{Q_1'} \dots J_{Q_N'}}],$$
	where now 
	$$J_Q = \int_{Q\times Q} :e^{i \beta \hat{Y}_\delta(x)}: :e^{i \beta \hat{Y}_\delta(y)}: R_\delta(x,y) dx dy.$$
	By Hölder's inequality we can bound 
	\[\E[J_{Q_1} \dots J_{Q_N} \overline{J_{Q_1'} \dots J_{Q_N'}}] \leq \E |J_{Q_1}|^{2N}.\]
	By scaling the right hand side equals
	\begin{align*}
	    & \delta^{4Nd} \int_{[0,1]^{4Nd}} \prod_{j=1}^N R_\delta(\delta x_j, \delta y_j) R_\delta(\delta x_j', \delta y_j') e^{\beta^2 \mathcal{E}(Y^{(\delta)};\mathbf{x},\mathbf{y};\mathbf{x}',\mathbf{y}')} \\
	    & \le \delta^{4Nd} \int_{[0,1]^{4Nd}} \prod_{j=1}^N \sqrt{\log \tfrac{C}{|x_j-\pi(x_j)|}\log \tfrac{C}{|y_j-\pi(y_j)|}\log \tfrac{C}{|x_j'-\pi(x_j')|}\log \tfrac{C}{|y_j'-\pi(y_j')|}} e^{\beta^2 \mathcal{E}(Y^{(\delta)};\mathbf{x},\mathbf{y};\mathbf{x}',\mathbf{y}')},
	\end{align*}
	where we have used Lemma~\ref{lemma:star_scale_cov_2} and $\pi(x)$ denotes the closest point to point $x$ in the set \[\{x_1,\dots,x_N,y_1,\dots,y_N,x_1',\dots,x_N',y_1',\dots,y_N'\}\setminus\{x\}.\]
	By relabeling the points as $z_1,\dots,z_{4N}$ and using Lemma~\ref{lemma:star_scale_onsager} we then have the upper bound
	\[\delta^{4Nd} \int_{[0,1]^{4Nd}} \prod_{j=1}^{4N} \sqrt{\log \frac{C}{|z_j - z_{F(j)}|}} \frac{1}{|z_j - z_{F(j)}|^{\beta^2/2}},\]
	which by Lemma~\ref{lem:integral_min} is bounded by
	\[C^N (d-\beta^2)^{-4N} \delta^{4Nd} N^{4N}\]
	for some constant $C > 0$.
	Hence we can bound $\E|\scalar{\overline{DM}, h_\delta}|^{2N} $ by
	\[C^N(d-\beta^2)^{-4N}\delta^{4Nd}N^{4N}\beta^{2N}\delta^{2\beta^2N}e^{2c\sqrt{\delta} N^2} \delta^{-2Nd}\int_{ K^{2N}} \exp \left( 4 \beta^2 \mathcal{E}( Y_{\delta^{1/2}}; \mathbf{x} ; \mathbf{x}' )\right), \]
	where for convenience we have turned $\mathbf{q}, \mathbf{q}'$ back to $\mathbf{x}, \mathbf{x}'$ by paying the same price.
	The latter integral is the $2N$-th moment of the $2\beta$ chaos of field $Y_{\delta^{1/2}}$, which by Lemma~\ref{lemma:star_scale_onsager} and \eqref{eq:super_critical_blowup} is bounded by $C^N N^{2N}\Big(\log \frac{1}{\delta}\Big)^{N}\delta^{-N\max(2\beta^2 - \frac{d}{2},0)}$, giving
	\[\E|\scalar{\overline{DM}, h_\delta}|^{2N} \leq C^Ne^{c\sqrt{\delta} N^2}(d-\beta^2)^{-4N}\Big(\log \frac{1}{\delta}\Big)^{N}\delta^{N(2d+\min(\frac{d}{2},2\beta^2))}N^{6N}.\]
	Note that for any fixed $b,C,\nu>0$ we have $2 b^{-1} C \log \frac{1}{\delta} < \delta^{-\nu}$ and $\delta$ small enough.
	One thus sees that
	\[\P[|\langle \overline{DM}, h_\delta \rangle_H| \ge b (d-\beta^2)^{-2}\delta^{d}] \le 2^{-N} e^{c \sqrt{\delta} N^2} \delta^{-\nu N} \delta^{N \min(\frac{d}{2}, 2\beta^2)} N^{6N}\]
	yields the desired upper bound by choosing e.g. $N = \delta^{-\beta^2/(24d)}$.
	\end{proof}

\subsubsection{Special case: the standard log-correlated field on the circle}\label{sec:circlecase}

In this section we will briefly explain how to extend the proof of Proposition~\ref{prop:tail_det} to the case where we are interested in the total mass of the imaginary chaos defined using the field $\Gamma$ on the unit circle which has the covariance $\log \frac{1}{|x-y|}$, where one now thinks of $x$ and $y$ as being complex numbers of modulus $1$. See Section \ref{sec:basic} for the precise definitions. 

Recall, that the extra complication in this case is that the field is degenerate in the sense that it is conditioned to satisfy $\int_0^1 \Gamma(e^{2 \pi i \theta}) \, d\theta = 0$. In terms of the proof of Proposition~\ref{prop:tail_det} this creates some annoyance, as the function $h_\delta$ we used in the projection bounds does not anymore belong to the Cameron--Martin space $H_\C$ of $\Gamma$, and we will instead need to look at the function $\tilde{h}_\delta = h_\delta - \int h_\delta(y) \, dy$.

As the field $\Gamma(e^{2\pi i \cdot})$ is non-degenerate when restricted to $I_0 \coloneqq [-1/4,1/4]$ (see again Section \ref{sec:basic}), it is also beneficial to introduce a smooth bump function $\psi$ supported in $I_0 \coloneqq [-1/4,1/4]$ , and thus set
\[h_\delta(x) = \psi(x) e^{i\beta Y_\delta(x) - \frac{\beta^2}{2} \E[Y_\delta(x)^2]} \int_{I_0} \psi(y) :e^{i\beta (\hat{Y}_\delta(y)+Z(y))}: R_\delta(x,y) \, dy.\]
This will let us still use the decomposition $X = Y + Z$ where $\Gamma_{|I_0} = X_{|I_0}$ and streamline most of the proof.

In the case of Lemmas \ref{lemma:mainterm} and \ref{lemma:oddterm}, i.e. in terms $\langle DM, \tilde{h}_\delta \rangle_{H_\C}$ and $\langle \overline{DM}, \tilde{h}_\delta \rangle_{H_\C}$, this subtraction of the mean introduces the extra term $i\beta M \int_0^1 h_\delta(y) \, dy$. In the case of Lemma \ref{lemma:normterm}, we have an extra term of the form $|\int_0^1 h_\delta(y)|$. The next lemma guarantees that both terms are negligible.

\begin{lemma}
For all $c > 0$ there is some $c_1 > 0$ such that we have 
$$\P[|\int_0^1 h_\delta(y) \, dy| > c \delta (1-\beta^2)^{-1/2}] \le e^{-c_1\delta^{-1} c^{\frac{2}{\beta^2}}}$$ and $$\P[|M \int_0^1 h_\delta(y) \, dy| > c \delta (1-\beta^2)^{-1}] \le e^{-c_1\delta^{-1/2} c^{\frac{1}{\beta^2}}}$$
for all $\delta$ small enough.
\end{lemma}

\begin{proof}
   We will bound the $N$--th moment of $|M \int h_\delta(y)|$, use the Chebyshev inequality and optimize over $N$.
   Note that by the Cauchy--Schwarz inequality we have
   \[\Expect{ \abs{ M \int_0^1 h_\delta(y) \, dy}^N } \le \E[|M|^{2N}]^{1/2} \Expect{ \abs{ \int_0^1 h_\delta(y) \, dy}^{2N}}^{1/2}\]
   and by \cite[Theorem~1.3]{JSW} we know that (recall that we are currently in a one-dimensional setting)
   \[\E[|M|^{2N}] \le C^N (d-\beta^2)^{-N}N^{\beta^2 N}\]
   for some $C > 0$. We mention that, in the article \cite{JSW}, the dependence of the above constant in terms of $\beta$ was not stated but follows from their approach (see \eqref{eq:lem_bound_integral_1}).
   To bound $\E[|\int_0^1 h_\delta(y) \, dy|^{2N}]$, we note that by Jensen's inequality we have
   \[\E\Big[\Big|\int_0^1 h_{\delta}(y) \, dy\Big|^{2N}\Big] \le \E\Big[\Big(\int_0^1 |h_{\delta}(y)|^2 \, dy\Big)^N\Big],\]
   where the right hand side equals
\[\E\Big[\Big(\int_0^1 |\psi(x)|^2 e^{-\beta^2 \E[Y_\delta(x)^2]} \Big|\int_0^1 \psi(y) :e^{i\beta (\hat{Y}_\delta(y) + Z(y))}: R_\delta(x,y) \, dy\Big|^2 \, dx\Big)^N\Big]. \]
We bound $|\psi(x)|^2 e^{-\beta^2 \E[Y_\delta(x)^2]}$ by $C \delta^{\beta^2}$ and since $R_\delta(x,y) = 0$ whenever $x,y$ do not belong to the same square, we can bound the above expression by 
\[C^N\delta^{N\beta^2}\delta^{-N} \sum_{Q \in \Qc_\delta} \E \Big[\Big(\int_{Q^3}\psi(y)\psi(z) :e^{i\beta(\hat{Y}_\delta(y)+Z(y))}:  R_\delta(x,y)R_\delta(x,z)  :e^{-i\beta(\hat{Y}_\delta(z)+Z(z))}:\, dz \, dx \, dy \Big)^N\Big].\]
By developing the expectation into a multiple integral, using an Onsager inequality associated to the smooth field $Z$ (see \eqref{eq:onsager_smooth}) and then rewriting the multiple integrals as an expectation, we see that we can get rid of the field $Z$ in the above expectation by only paying a multiplicative price $C^N$. 

Thus it remains to bound
\[
C^N\delta^{N\beta^2}\delta^{-N} \sum_{Q \in \Qc_\delta} \E \Big[\Big(\int_{Q^3}\psi(y)\psi(z) :e^{i\beta\hat{Y}_\delta(y)}:  R_\delta(x,y)R_\delta(x,z)  :e^{-i\beta\hat{Y}_\delta(z)}:\, dz \, dx \, dy \Big)^N\Big].
\]
By scaling we see that each term in the sum is equal in law to
\[\delta^{3N} J_Q \coloneqq \delta^{3N}  \E \Big[\Big(\int_{\delta^{-1}Q \times \delta^{-1}Q  \times \delta^{-1}Q} \psi(\delta y) \psi(\delta z):e^{i\beta \hat{Y}_\delta(\delta y)}:  R_\delta(\delta x, \delta y)R_\delta(\delta x,\delta z)  :e^{-i\beta \hat{Y}_\delta(\delta z)}: \, dz \, dx \, dy \Big)^N\Big].\]
To bound this expectation, we expand the product and obtain a multiple integral over $x_i, y_i, z_i$, $i=1 \dots N$. The expectation of the product of $:e^{i\beta \hat{Y}_\delta(\delta y)}:$ and $:e^{-i\beta \hat{Y}_\delta(\delta z)}:$ leads to $\mathcal{E}(\hat{Y}_\delta(\delta \cdot) ; \mathbf{y} ; \mathbf{z})$ that we bound using the Onsager inequality \eqref{eq:onsager_tail}. Since for any fixed $y$ and $z$,
\[
\psi(\delta y) \psi( \delta z) \int_{\delta^{-1}Q} R_\delta(\delta x, \delta y) R_\delta(\delta x, \delta z)\, dx < C,
\]
we can first integrate the variables $x_i$ and control the remaining integral over $y_i$ and $z_i$, $i=1 \dots N$ with \eqref{eq:lem_bound_integral_1}. Overall, $J_Q$ is bounded by $(d-\beta^2)^{-N}N^{\beta^2N}$.

Altogether we obtain that
 \[\E \Big[\Big|\int_0^1 h_\delta(y) \, dy\Big|^{2N}\Big] \le C^N (d-\beta^2)^{-N}\delta^{(\beta^2+2)N}N^{\beta^2 N}\]
and hence
   \[\E \Big[\Big|M \int_0^1 h_\delta(y) \, dy\Big|^N \Big] \le C^N (d-\beta^2)^{-N} \delta^{(\frac{\beta^2}{2} + 1)N} N^{\beta^2 N},\]
   which gives us the tail estimates
      \[\P \Big[ \Big|\int_0^1 h_\delta(y) \, dy \Big| \ge \lambda (d-\beta^2)^{-1/2} \Big] \le \frac{C^N \delta^{(\frac{\beta^2}{2} + 1)N} N^{\frac{\beta^2}{2}N}}{\lambda^N}.\]
and
   \[\P \Big[ \Big| M \int_0^1 h_\delta(y) \, dy \Big| \ge \lambda (d-\beta^2)^{-1} \Big] \le \frac{C^N \delta^{(\frac{\beta^2}{2} + 1)N} N^{\beta^2 N}}{\lambda^N}.\]
   Optimising over $N$ now concludes.
\end{proof}

\appendix

\section{Appendix: some standard proofs}\label{appendixA}

\begin{proof}[Proof of Lemma \ref{lemma:star_scale_cov_2}]
It is calculationally somewhat easier to work with the  rescaled field $Y^{(\epsilon)}(x) = \hat{Y}_\epsilon(\delta x)$, which
can be expressed using white noise as:
\[Y^{(\delta)}(x) \coloneqq \int_{\reals^d \times [0,\infty)} e^{du/2} \tilde{k}(e^u(t-x))\sqrt{1 - \delta^{\alpha} e^{-\alpha u}} dW(t,u).\]
The first inequality then follows directly:
  \[\E[Y^{(\delta)}(x)Y^{(\delta)}(y)] = \int_0^\infty k(e^u(x-y))(1 - \delta^{\alpha}e^{-\alpha u}) \, du \le \int_0^\infty k(e^u(x-y)) \, du \le \log \frac{1}{|x-y|}\]
  by the fact that $k$ is supported in $B(0,1)$ and $k(t) \le 1$ for all $t$.

  For the second inequality we compute
  \begin{align*}
    \int_0^\infty k(e^u(x-y))(1 - \delta^\alpha e^{-\alpha u}) \, du & \ge \int_0^\infty k(e^u(x-y))(1 - e^{-\alpha u}) \, du \\
    & \ge \int_0^{\log \frac{1}{|x-y|}} k(e^u(x-y)) \, du - \int_0^{\infty} e^{-\alpha u} \, du \\
    & \ge \log \frac{1}{|x-y|} + \int_0^{\log \frac{1}{|x-y|}} (k(e^u(x-y)) - 1) \, du - \frac{1}{\alpha}
  \end{align*}
  Note that by Taylor's theorem we have for all $t \in \reals$ the inequality
  \[k(t) \ge 1 + k'(0) t - c t^2\]
  for some constant $c > 0$, and in fact since $k$ is smooth and symmetric we have $k'(0) = 0$.
  Hence
  \[\int_0^{\log \frac{1}{|x-y|}} (k(e^u(x-y)) - 1) \, du \ge -c \int_0^{\log \frac{1}{|x-y|}} e^{2u}|x-y|^2 = -c (\frac{1}{2|x-y|^2} |x-y|^2 - \frac{|x-y|^2}{2}) \ge -\frac{c}{2},\]
  from which the claim follows.
  
  Finally, the independence comes from the fact that $k$ is supported in $B(0,1)$
\end{proof}

\begin{proof}[Proof of Lemma~\ref{lemma:star_scale_onsager}]
Let us begin with the field $Y_\varepsilon$.
Set $q_j = 1$ for $1 \le j \le N$ and $q_j = -1$ for $N+1 \le j \le 2N$ and note that
\[\mathcal{E}(Y_\varepsilon;\mathbf{x};\mathbf{y}) = -\frac{1}{2} \E\left[\Big(\sum_{j=1}^{2N} q_j Y_{d_j \wedge \varepsilon}(z_j)\Big)^2\right] + \frac{1}{2}\sum_{j=1}^{2N} \E[Y_{d_j \wedge \varepsilon}(z_j)^2] \le \frac{1}{2}\sum_{j=1}^{2N} \log \frac{1}{d_j \wedge \varepsilon}\]
since $\E[Y_\varepsilon(x)Y_\varepsilon(y)] = \E[Y_s(x) Y_t(y)]$ for all $s,t \le \varepsilon \wedge |x-y|$ and $\E[Y_{\delta}(x)^2] \le \log \frac{1}{\delta}$ for all $\delta \in (0,1)$.

As the field $\hat{Y}_\varepsilon(\varepsilon x)$ has the same distribution as the field $Y^{(\varepsilon)}(x)$ from the proof of Lemma~\ref{lemma:star_scale_cov_2}, we have
\[\mathcal{E}(\hat{Y}_\varepsilon(\varepsilon \cdot);\mathbf{x};\mathbf{y}) = -\frac{1}{2} \E\left[\Big(\sum_{j=1}^{2N} q_j \hat{Y}^{(\varepsilon)}_{d_j}(z_j)\Big)^2\right] + \frac{1}{2} \sum_{j=1}^{2N} \E[\hat{Y}^{(\varepsilon)}_{d_j}(z_j)^2] \le \frac{1}{2} \sum_{j=1}^{2N} \log \frac{1}{d_j}.\]

Finally, if $R$ is a regular field then
\[\mathcal{E}(R;\mathbf{x};\mathbf{y}) = -\frac{1}{2} \E\left[\Big(\sum_{j=1}^{2N} q_j R(z_j)\Big)^2\right] + \frac{1}{2} \sum_{j=1}^{2N} \E[R(z_j)^2] \le N \sup_{1 \le j \le 2N} \E[R(z_j)^2].\qedhere\]
\end{proof}

\begin{proof}[Proof of Lemma \ref{lem:Nualart}]
We prove this lemma in the context of real-valued random variables. The extension to complex-valued random variables follows immediately.

In page 58 of \cite{Nualart}, an operator $L$ on the set of variables with finite second moment is introduced and used to define the norm $\norme{| F |}_{k,p} := \Expect{ ((I-L)^{k/2} F)^p}^{1/p}$.
The norms $\norme{| \cdot |}_{k,p}$ and $\norme{\cdot}_{k,p}$ are equivalent (see \cite{Nualart} page 77). Hence $\sup_n \Expect{ ((I-L)^{k/2} F_n)^p} < \infty$. By weak compactness of balls in $L^p(\Omega)$, we can extract a subsequence $(n(i), i \geq 1)$ such that $((I-L)^{k/2} F_{n(i)}, i \geq 1)$ converges weakly towards some element $G$. Since the $L^p$-norm is weakly lower-semicontinuous, we moreover have
\[
\Expect{G^p} \leq \liminf_i \Expect{((I-L)^{k/2} F_{n(i)})^p} \leq \limsup_n \Expect{((I-L)^{k/2} F_n)^p}.
\]
In the proof of \cite[Lemma 1.5.3]{Nualart}, D. Nualart shows that $F = (I-L)^{-k/2} G$. This implies that
\[
\norme{F}_{k,p} \leq C_{k,p} \norme{|F|}_{k,p} = C_{k,p} \Expect{G^p}^{1/p} \leq C_{k,p} \limsup_n \norme{|F_n|}_{k,p} \leq C_{k,p}' \limsup_n \norme{F_n}_{k,p}.
\]
This concludes the proof.
\end{proof}

\section{Appendix: proof of Proposition~\ref{prop:dinfty}}\label{appendixB}

\begin{proof}[Proof of Proposition~\ref{prop:dinfty}]
We start by showing that $M$ belongs to $\D^\infty$. Let $n \geq 1, \delta > 0, j \geq 0$ and $p \geq 1$.
In the following, we will denote
\[
\Gamma_\delta = \Gamma * \varphi_\delta, \quad \Gamma_{n,\delta} = \sum_{k=1}^n A_k e_k \ast \varphi_\delta,
\quad
M_\delta = \int_\C f(x) e^{i \beta \Gamma_\delta(x) + \frac{\beta^2}{2} \E[\Gamma_\delta(x)^2]} dx
\]
and
\[
M_{n,\delta} = \int_\C f(x) e^{i \beta \Gamma_{n,\delta}(x) + \frac{\beta^2}{2} \E[\Gamma_{n,\delta}(x)^2] } dx.
\]
$M_{n,\delta}$ is a smooth random variable and $D^j M_{n,\delta}$ is equal to
\begin{equation}
(i \beta)^j \int_\C dx f(x) e^{i \beta \Gamma_{n,\delta}(x) + \frac{\beta^2}{2} \E[\Gamma_{n,\delta}(x)^2]} \sum_{k_1, \dots, k_j =1}^n (e_{k_1} \ast \varphi_\delta)(x) \dots (e_{k_j} \ast \varphi_\delta)(x) e_{k_1} \otimes \dots \otimes e_{k_j}. \label{eq:proof_lem_density0}
\end{equation}
Since $(e_{k_1} \otimes \dots \otimes e_{k_j}, k_1, \dots, k_j = 1 \dots n)$ is an orthonormal family of $H^{\otimes j}$, we deduce that
\begin{align*}
& \norme{D^j M_{n,\delta} }_{H_\C^{\otimes j}}^2
= \beta^{2j} \int_{\C^2} f(x) f(y) e^{i \beta \Gamma_{n,\delta}(x) - i\beta \Gamma_{n,\delta}(y) + \frac{\beta^2}{2} \E [\Gamma_{n,\delta}(x)^2] + \frac{\beta^2}{2} \E [\Gamma_{n,\delta}(y)^2]} \\
& \quad \times \left( \sum_{k=1}^n (e_k \ast \varphi_\delta)(x) (e_k \ast \varphi_\delta)(y) \right)^j dx dy.
\end{align*}
Thanks to the convolution, all the integrated terms are uniformly bounded in $n$ and $x_1 \dots x_{p}$, $y_1 \dots y_{p}$. By dominated convergence theorem and then by using \eqref{eq:lem_Onsager_conv} which provides an Onsager inequality for convolution approximations, we deduce that
\begin{align*}
& \limsup_{n \to \infty} \Expect{ \norme{D^j M_{n,\delta} }_{H_\C^{\otimes j}}^{2p} } \\
& \leq \beta^{2jp} \int_{\C^{2p}} dx_1 \dots dx_{p} dy_1 \dots dy_{p}
\prod_{l=1}^{p} f(x_l) f(y_l) \left( C \ast (\varphi_\delta \otimes \varphi_\delta)(x_l,y_l) \right)^j e^{\beta^2 \Ec(\Gamma_\delta;\mathbf{x}; \mathbf{y})} \\
& \leq C_{j,p} \norme{f}_\infty^{2p} \int_{K^{2p}} dz_1 \dots dz_{2p} \prod_{l=1}^{2p} \left( \min_{l' \neq l} \abs{z_l - z_{l'}} \right)^{-\beta^2/2} \left(\max_{l' \neq l} C \ast (\varphi_\delta \otimes \varphi_\delta)(z_l,z_{l'}) \right)^{j/2}
\end{align*}
where $K$ is the support of $f$. Importantly, the above constant $C_{j,p}$ does not depend on $\delta$. Notice that
\[
C \ast (\varphi_\delta \otimes \varphi_\delta)(x,y) \leq C \log \frac{c}{|x-y| \vee \delta}.
\]
Hence, if we let $\eps >0$ be such that $\beta^2/2 + \eps < d/2$, there exists $C_{j,p}'>0$ independent of $\delta$ such that
\begin{align}
\limsup_{n \to \infty} \Expect{ \norme{D^j M_{n,\delta} }_{H_\C^{\otimes j}}^{2p} } \leq C_{j,p}' \int_{K^{2p}} dz_1 \dots dz_{2p} \prod_{l=1}^{2p} \left( \min_{l' \neq l} \abs{z_l - z_{l'}} \right)^{-\beta^2/2-\eps} \leq C_{j,p}'' \label{eq:proof_lem_density}
\end{align}
by \eqref{eq:lem_bound_integral_1}.
Since $(M_{n,\delta}, n \geq 1)$ converges in $L^{2p}$ towards $M_\delta$,
Lemma \ref{lem:Nualart} and \eqref{eq:proof_lem_density} imply that for all $k \geq 1$, $M_\delta \in \D^{k,2p}$ and that 
\begin{equation}
\label{eq:proof_lem_density8}
\sup_{\delta > 0} \norme{M_\delta}_{k,2p} < \infty.
\end{equation}
Now, because $(M_\delta, \delta >0)$ converges in $L^{2p}$ towards $M$, Lemma \ref{lem:Nualart} implies that for all $k \geq 1$, $M \in \D^{k,2p}$. This concludes the proof that $M \in \D^\infty$.\\

We now turn to the proof of the formula for $DM$. On the one hand, \eqref{eq:proof_lem_density0} gives
\[
DM_{n,\delta} = i \beta \int_\C dx f(x) e^{i \beta \Gamma_{n,\delta}(x) + \frac{\beta^2}{2} \E [\Gamma_{n,\delta}(x)^2]} \sum_{k=1}^n (e_k \ast \varphi_\delta)(x) e_k.
\]
One can then show that $(DM_{n,\delta},n \geq 1)$ converges in $L^2(\Omega;H)$ towards
\[
i \beta \int_\C dx f(x) e^{i\beta \Gamma_\delta(x) + \frac{\beta^2}{2} \E [\Gamma_\delta(x)^2]} \sum_{k=1}^\infty (e_k \ast \varphi_\delta)(x) e_k.
\]
On the other hand, the first part of the proof showed that $\sup_n \Expect{ \norme{ DM_{n,\delta} }_{H_\C}^2 } < \infty$ and Lemma \ref{lem:Nualart0} implies that $(DM_{n,\delta}, n \geq 1)$ converges to $DM_\delta$ in the weak topology of $L^2(\Omega;H)$. Hence
\[
DM_\delta = i \beta \int_\C dx f(x) e^{i\beta \Gamma_\delta(x) + \frac{\beta^2}{2} \E[\Gamma_\delta(x)^2]} \sum_{k=1}^\infty (e_k \ast \varphi_\delta)(x) e_k.
\]
Let us now show that $(DM_\delta, \delta>0)$ converges in $L^2(\Omega;H)$ towards
\[
i \beta  \int_\C dx f(x) \mu(x) C(x, \cdot).
\]
Firstly, since
\[
C(x,\cdot) = \sum_{k \geq 1} e_k(x) e_k(\cdot)
\]
and the $e_k, k \geq 1,$ form an orthonormal family of $H$, we have
\begin{align}
& \Expect{ \norme{ \int_\C dx f(x) \mu(x) C(x, \cdot) - \int_\C dx f(x) e^{i \beta \Gamma_\delta(x) + \frac{\beta^2}{2} \E[\Gamma_\delta(x)^2]} C(x,\cdot) }_{H_\C}^2 } \label{eq:proof_lem_density2} \\
& = \sum_{k \geq 1} \Expect{ \left( \int_\C f(x) \mu(x) e_k(x) dx - \int_\C f(x) e^{i \beta \Gamma_\delta(x) + \frac{\beta^2}{2} \E[\Gamma_\delta(x)^2]} e_k(x) dx \right)^2 }. \nonumber
\end{align}
Each single term in the above sum goes to zero as $\delta \to 0$. Moreover, using Onsager inequality for convolution approximations \eqref{eq:lem_Onsager_conv}, one can obtain a domination in a similar manner as what we did in the first part of the proof. By the dominated convergence theorem, it implies that \eqref{eq:proof_lem_density2} goes to zero as $\delta \to 0$.
Secondly,
\begin{align}
& \Expect{ \norme{ \int_\C dx f(x) e^{i \beta \Gamma_\delta(x) + \frac{\beta^2}{2} \E[\Gamma_\delta(x)^2]} \sum_{k \geq 1} (e_k \ast \varphi_\delta)(x) e_k - \int_\C dx f(x) e^{i \beta \Gamma_\delta(x) + \frac{\beta^2}{2} \E[\Gamma_\delta(x)^2]} C(x,\cdot)  }_{H_\C}^2 } \label{eq:proof_lem_density3} \\
& = \sum_{k \geq 1} \Expect{ \left( \int_\C f(x) e^{i \beta \Gamma_\delta(x) + \frac{\beta^2}{2} \E[\Gamma_\delta(x)^2]} ((e_k \ast \varphi_\delta)(x) - e_k(x)) dx \right)^2} \nonumber \\
& \leq C \norme{f}_\infty^2 \int_{K^2} \abs{x-y}^{-\beta^2} \abs{ \sum_{k \geq 1} ((e_k \ast \varphi_\delta)(x) - e_k(x))((e_k \ast \varphi_\delta)(y) - e_k(y)) } dx dy \nonumber
\end{align}
where $K$ is as before the support of $f$.
The above integrand is dominated by the integrable function $C \abs{x-y}^{-\beta^2} \log(c/|x-y|)$. Dominated convergence theorem thus implies that \eqref{eq:proof_lem_density3} goes to zero as $\delta \to 0$. Putting things together, we have shown the aforementioned convergence: $(DM_\delta, \delta>0)$ converges in $L^2(\Omega;H)$ towards
\[
i \beta \int_\C dx f(x) \mu(x) C(x, \cdot).
\]
With \eqref{eq:proof_lem_density8}, we notice that $\sup_\delta \Expect{ \norme{DM_\delta}_{H_\C}^2 } < \infty$ and Lemma \ref{lem:Nualart0} also shows that $(DM_\delta, \delta >0)$ converges to $DM$ in the weak topology of $L^2(\Omega;H)$. This yields
\[
DM = i \beta \int_\C dx f(x) \mu(x) C(x, \cdot).\qedhere
\]
\end{proof}


\begin{thebibliography}{99}
  \bibitem{AdamsFrazier} Adams, D.R., Frazier, M.: Composition operators on potential spaces. P. Am. Math. Soc. 114(1), 155--165 (1992). \url{https://doi.org/10.2307/2159794}
  \bibitem{AdlerTaylor} Adler, R.J., Taylor, J.E.: Random fields and geometry. Springer, New York (2007). \url{https://doi.org/10.1007/978-0-387-48116-6}
  \bibitem{Aronz} Aronszajn, N.: Theory of reproducing kernels. T. Am. Math. Soc. 68(3), 337--404 (1950). \url{https://doi.org/10.1090/S0002-9947-1950-0051437-7}
  \bibitem{AruJunnila} Aru, J., Junnila, J.: Reconstructing the base field from imaginary multiplicative chaos. B. Lond. Math. Soc. (2021). \url{https://doi.org/10.1112/blms.12466}
  \bibitem{BJM} Barral, J., Jin, X., Mandelbrot, B.: Convergence of complex multiplicative cascades. Ann. Appl. Probab. 20(4), 1219--1252 (2010). \url{https://doi.org/10.1214/09-AAP665}
  \bibitem{BM} Barral, J., Mandelbrot, B.: Fractional multiplicative processes. Ann. I. H. Poincaré B. 45(4), 1116--1129 (2009). \url{https://doi.org/10.1214/08-AIHP198}
  \bibitem{BH} Behzadan, A., Holst, M.: Multiplication in Sobolev spaces, revisited. arXiv:1512.07379
  \bibitem{Big} Biggins, J.D.: Uniform Convergence of Martingales in the Branching Random Walk. Ann. Probab. 20(1), 137--151 (1992). \url{https://doi.org/10.1214/aop/1176989921}
  \bibitem{BogachevGauss} Bogachev, V.I.: Gaussian measures. Mathematical Surveys and Monographs, 62. American Mathematical Society, Providence, RI (1998). \url{https://doi.org/10.1090/surv/062}
  \bibitem{CGPR} Camia, F., Gandolfi, A., Peccati, G., Reddy, T.R.: Brownian Loops, Layering Fields and Imaginary Gaussian Multiplicative Chaos. Commun. Math. Phys. 381(3), 889--945 (2021). \url{https://doi.org/https://doi.org/10.1007/s00220-020-03932-9}
  \bibitem{ChhNaj} Chhaibi, R., Najnudel, J.: On the circle, $GMC^\gamma = \varprojlim C\beta E_n$ for $\gamma = \sqrt{\frac{2}{\beta}}, $ $( \gamma \leq 1 )$. arXiv: 1904.00578.
  \bibitem{DERR} Derrida, B., Evans, M.R., Speer, E.R.: Mean field theory of directed polymers with random complex weights. Commun. Math. Phys. 156(2), 221-244 (1993). \url{https://doi.org/10.1007/BF02098482}
  \bibitem{FBoriginal} Fyodorov, Y.V., Bouchaud, J.--P.: Freezing and extreme-value statistics in a random energy model with logarithmically correlated potential. J. Phys. A--Math. Theor. 41(37), 372001 (2008). \url{https://doi.org/10.1088/1751-8113/41/37/372001}
  \bibitem{GarbanSepulveda} Garban, C., Sepúlveda, A.: Statistical reconstruction of the Gaussian free field and KT transition. arXiv:2002.12284.
  \bibitem{JSW2} Junnila, J., Saksman, E., Webb, C.: Decompositions of log-correlated fields with applications. Ann. Appl. Probab. 29(6), 3786--3820 (2019). \url{https://doi.org/10.1214/19-AAP1492}
  \bibitem{JSW} Junnila, J., Saksman, E., Webb, C.: Imaginary multiplicative chaos: Moments, regularity and connections to the Ising model. Ann. Appl. Probab. 30(5), 2099--2164 (2020). \url{https://doi.org/10.1214/19-AAP1553}
  \bibitem{DOZZ} Kupiainen, A., Rhodes, R., Vargas, V.: Integrability of Liouville theory: proof of the DOZZ formula. Ann. Math. 191(1), 81--166 (2020). \url{https://doi.org/10.4007/annals.2020.191.1.2}
  \bibitem{LCRV} Lacoin, H., Rhodes, R., Vargas, V.: Complex Gaussian multiplicative chaos. Commun. Math. Phys. 337(2), 569--632 (2015). \url{https://doi.org/10.1007/s00220-015-2362-4} 
  \bibitem{LRV} Lacoin, H., Rhodes, R., Vargas, V.: A probabilistic approach of ultraviolet renormalisation in the boundary Sine-Gordon model. arXiv:1903.01394
  \bibitem{LSZ} Leblé, T., Serfaty, S., Zeitouni, O.: Large deviations for the two-dimensional two-component plasma. Commun. Math. Phys. 350(1) 301--360 (2017). \url{https://doi.org/10.1007/s00220-016-2735-3}
  \bibitem{Malliavin} Malliavin, P.: Stochastic calculus of variations and hypoelliptic operators. Proc. Internat. Symposium on Stoch. Diff. Equations (1976), Kyoto Univ. Press, Wiley, 195--263 (1978).
  \bibitem{Nualart} Nualart, D.: The Malliavin calculus and related topics. Springer, Berlin (2006). \url{https://doi.org/10.1007/3-540-28329-3}
  \bibitem{NuaNua} Nualart, D., Nualart, E.: Introduction to Malliavin Calculus. Cambridge University Press, Cambridge (2018). \url{https://doi.org/10.1017/9781139856485}
  \bibitem{EPreview} Powell, E.: Critical Gaussian multiplicative chaos: a review. arXiv: 2006.13767
  \bibitem{FB} Remy, G.: The Fyodorov--Bouchaud formula and Liouville conformal field theory. Duke Math. J. 169(1), 177--211 (2020). \url{https://doi.org/10.1215/00127094-2019-0045}
  \bibitem{RVreview} Rhodes, R., Vargas, V.: Gaussian multiplicative chaos and applications: A review. Probab. Surveys 11, 315--392 (2014). \url{https://doi.org/10.1214/13-PS218}
  \bibitem{RV} Robert, R., Vargas, V.: Gaussian multiplicative chaos revisited. Ann. Probab.  38(2), 605--631 (2010). \url{https://doi.org/10.1214/09-AOP490}
  \bibitem{SW} Saksman, E., Webb, C.: The Riemann zeta function and Gaussian multiplicative chaos: Statistics on the critical line. Ann. Probab. 48(6), 2680--2754 (2020). \url{https://doi.org/10.1214/20-AOP1433}
  \bibitem{SSV} Schoug, L., Sep\'ulveda, A., Viklund, F.: Dimension of two-valued sets via imaginary chaos. Int. Math. Res. Notices (2020). \url{https://doi.org/10.1093/imrn/rnaa250}
  \bibitem{Triebel} Triebel, H.: Interpolation Theory, Function Spaces, Differential Operators. North-Holland Mathematical Library 18. North-Holland, Amsterdam (1978).
  \bibitem{VanHandel} van Handel, R.: Probability in High Dimension. APC550 Lecture notes.
\end{thebibliography}
\end{document}